\documentclass[12pt]{article}
\usepackage{amsmath,amssymb, amsthm, graphicx}

\def\ga{c}

\newtheorem{thm}{Theorem}
\newtheorem{conj}[thm]{Conjecture}
\newtheorem{prop}[thm]{Proposition}
\newtheorem{lem}[thm]{Lemma}
\newtheorem{lemma}[thm]{Lemma}

\newcommand{\id}{\operatorname{Id}}

\newcommand{\dd}[1]{{\bf #1}}

\title{More refined enumerations of alternating sign matrices}
\author{Ilse Fischer\thanks{Supported by the Austrian Science Foundation
    (FWF), grant number S9607--N13, in the framework of the National Research
    Network  ``Analytic Combinatorics and Probablistic Number Theory''.} \and Dan Romik\thanks{Supported by the Israel Science Foundation
(ISF) grant number 1051/08.}}
\begin{document}
\maketitle
\begin{center}
\end{center}
\abstract{We study a further refinement of the standard refined enumeration of alternating sign matrices (ASMs) according to their first two rows instead of just the first row, and more general ``$d$-refined'' enumerations of ASMs according to the first $d$ rows. For the doubly-refined case of $d=2$, we derive a system of linear equations satisfied by the doubly-refined enumeration numbers $A_{n,i,j}$ that enumerate such matrices. We give a conjectural explicit formula for $A_{n,i,j}$ and formulate several other conjectures about the sufficiency of the linear equations to determine the $A_{n,i,j}$'s and about an extension of the linear equations to the general
$d$-refined enumerations.
}

\section{Introduction}

\subsection{Alternating sign matrices}

An \dd{alternating sign matrix (ASM)} of order $n$ is an $n\times n$ matrix with entries in $\{0,-1,1\}$ such that in every row and every column, the sum of the entries is $1$ and the non-zero numbers appear with alternating signs; see Fig.\ 1(a) for an example with $n=5$. Alternating sign matrices were first defined and studied in the early 1980's by David Robbins and Howard Rumsey in connection with Charles Dodgson's condensation method for computing determinants. They discovered, and later proved \cite{robbinsrumsey}, that the $\lambda$-determinant, a natural generalization of the determinant that arises from the condensation method, has an expansion as a sum over all ASMs of order $n$, just as the ordinary determinant has an expansion as a sum over permutation matrices. 

Robbins and Rumsey were also interested in the \dd{enumeration} of ASMs. Together with William Mills, they denoted by $A_n$ the total number of ASMs of order $n$, and by $A_{n,k}$ for $1\le k\le n$  the number of ASMs with a $1$ in position $k$ of row 1, and conjectured \cite{millsetal2} two important enumeration identities,
nowadays known as the \dd{alternating sign matrix theorem} and the \dd{refined alternating sign matrix theorem}, which state that
\begin{eqnarray}
A_n &=& \prod_{j=0}^{n-1} \frac{(3j+1)!}{(n+j)!}\qquad\text{and} \label{eq:asmthm} \\
A_{n,k} &=& \binom{n+k-2}{k-1} \frac{(2n-k-1)!}{(n-k)!}\prod_{j=0}^{n-2} \frac{(3j+1)!}{(n+j)!}. \label{eq:refinedasmthm}
\end{eqnarray}
The attempt to prove these celebrated conjectures and subsequent developments have led to a rich combinatorial and algebraic theory relating ASMs to other combinatorial objects such as domino tilings and plane partitions, to representation theory and to the six-vertex model and other related lattice models in statistical physics. Eqs.\ \eqref{eq:asmthm} and \eqref{eq:refinedasmthm} were proved by Zeilberger \cite{zeilberger1, zeilberger2}, who in the refined case of \eqref{eq:refinedasmthm} 
(which contains \eqref{eq:asmthm} as a special case, since  $A_{n,1}=A_{n-1}$, and took longer to prove)
used the square-ice techniques introduced by Kuperberg \cite{kuperberg} in his simplified proof of \eqref{eq:asmthm}. 
An alternative proof of \eqref{eq:refinedasmthm} using entirely different methods
was found by Fischer \cite{fischer1, fischer2}.
Other variations on the Kuperberg-Zeilberger approach lead to different derivations of \eqref{eq:asmthm} by Colomo and Pronko \cite{colomopronkoasm} and of \eqref{eq:refinedasmthm} by Stroganov \cite{stroganov}.

Our goal in this paper is to extend these enumeration results to more refined enumerations of ASMs. In the simplest case, which we call the \dd{doubly-refined enumeration}, we will define numbers $A_{n,i,j}$ that enumerate ASMs based on their first \dd{two} rows rather than the first row\footnote{There is a different doubly-refined enumeration based on the first and \emph{last} rows, studied e.g.\ in \cite{colomopronkorefinedasm, fonsecaetal, stroganov}, which we will not discuss in this paper.}.
We will present a conjectural explicit formula for $A_{n,i,j}$ that turns out to be considerably more complicated than the product formulas \eqref{eq:asmthm} and \eqref{eq:refinedasmthm}. While proving this formula may look challenging, we take a significant first step by deriving and proving a system of linear equations satisfied (for each fixed $n$) by the numbers $(A_{n,i,j})_{i,j}$. Assuming a natural conjecture that we will formulate concerning the sufficiency of these equations to determine the $A_{n,i,j}$'s, proving the explicit formula is reduced to the (possibly quite complicated, but essentially mechanical) verification that the expression given by that formula satisfies the system of linear equations.

One can also consider the triply-refined, quadruply-refined, and more generally the \dd{$d$-refined enumeration} of ASMs based on their first $d$ rows. Our main result concerning the system of linear equations extends conjecturally (and empirically for $d=3$ and small values of $n$) to this generality. We believe our techniques are also relevant for attacking this more general problem (they are in fact an extension of the techniques used in \cite{fischer2} to attack the original refined case of $d=1$).
The existence of such detailed equations governing the $d$-refined enumeration is surprising and may offer a useful tool for attacking important open problems such as the problem of the limit shape of large random ASMs, for which a conjectural formula is given in \cite{colomopronko}.

\subsection{Monotone triangles}

To formulate our results precisely, we need first to define the doubly-refined enumeration numbers $A_{n,i,j}$. This is more conveniently done in terms of a family of objects called \dd{monotone triangles} that are related to ASMs. A monotone triangle of order $n$
is a triangular array
$ (t_{i,j})_{1\le i\le n, 1\le j\le i} $ of integers satisfying the inequalities
$$ t_{i,j} < t_{i,j+1}, \ \ \  t_{i,j} \le t_{i-1,j} \le t_{i,j+1} \qquad (2\le i\le n, 1\le j\le i-1). $$
In other words, in the coordinate system traditionally used to draw them, monotone triangles have strictly increasing rows
and weakly increasing northwest-southeast and southwest-northeast diagonals; see Fig.\ 1(b).

It is well-known (see \cite{bressoud, propp}) that ASMs correspond under a simple bijection to \dd{complete monotone triangles}, which are monotone triangles with the numbers $1,2,...,n$ in the bottom row. Formally, the bijection maps the ASM $M=(m_{i,j})$ to the monotone triangle
$T=(t_{i,j})$ such that $t_{i,1} < t_{i,2} < t_{i,3}< \ldots < t_{i,i}$ are the positions of the $1$'s in row $i$ of the \dd{column-sum matrix}
$ S=(s_{i,j})_{i,j=1}^n $
of $M$ defined by
$$ s_{i,j} = \sum_{r=1}^i m_{r,j}. $$
Using this bijection and the obvious symmetry with respect to vertical reflection in the definition of ASMs, it is easy to see that the enumeration and refined enumeration numbers $A_n$ and $A_{n,k}$ have the following equivalent definitions in terms of monotone 
triangles:
\begin{eqnarray}
A_n &=& \text{the number of complete monotone triangles of order $n$}, \nonumber \\
A_{n,k} &=& \text{the number of complete monotone triangles of order $n$ with the} \nonumber \\
& & \text{number $k$ in the first row} \nonumber \\
&=& \text{the number of monotone triangles of order $n-1$ with the} \nonumber \\
& & \text{numbers $1,2,\ldots,k-1,k+1,\ldots,n$ in the last row.} \label{eq:ankmonotone}
\end{eqnarray}

\begin{figure}
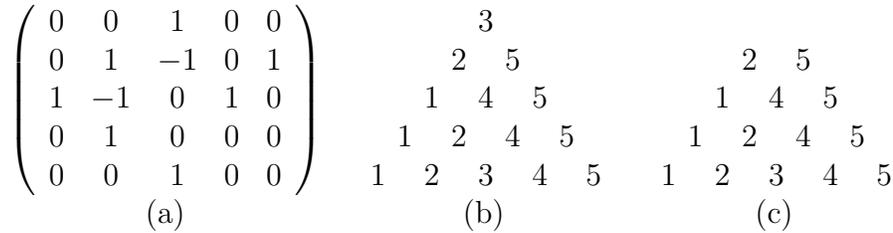

\begin{center}
\begin{tabular}{ccc}
$  \left( \begin{array}{ccccc}
0 & 0 & 1 & 0 & 0 \\
0 & 1 & -1 & 0 & 1 \\
1 & -1 & 0 & 1 & 0 \\
0 & 1 & 0 & 0 & 0 \\
0 & 0 & 1 & 0 & 0 
\end{array} \right)  
$
&
\begin{tabular}{*{13}{c@{\hspace{0.06in}}}}
 & && & 3 & & & &
 \\
& && 2 & & 5 & & &
 \\
&& 1 & & 4 & & 5 & &
 \\
 & 1 && 2 && 4 & & 5 &
 \\
 1 && 2 && 3 && 4 && 5
\end{tabular}
&
\begin{tabular}{*{13}{c@{\hspace{0.06in}}}}
 & && &  & & & &
 \\
& && 2 & & 5 & & &
 \\
&& 1 & & 4 & & 5 & &
 \\
 & 1 && 2 && 4 & & 5 &
 \\
 1 && 2 && 3 && 4 && 5
\end{tabular}  \\ 
(a)&(b)&(c)
\end{tabular}
\caption{An ASM of order 5, the corresponding complete monotone triangle, and the $(2,n)$-trapezoid obtained by deleting its first row.
\label{fig-illustr}}
\end{center}
\end{figure}

\subsection{The doubly-refined enumeration}

Our definition of $A_{n,i,j}$ is a natural generalization of \eqref{eq:ankmonotone}. For $n\ge 3$ and \hbox{$1\le i< j\le n$,} let
\begin{eqnarray*}
A_{n,i,j} &=& \text{the number of order $n-2$ monotone triangles with the}\\
& & \text{numbers $1,2,\ldots,\widehat{i},\ldots,\widehat{j},\ldots n$ in the last row,}
\end{eqnarray*}
where the notation $\widehat{i},\widehat{j}$ means that $i$ and $j$ are omitted from the list.
Equivalently, we can look at the second row of a monotone triangle and write
\begin{eqnarray*}
A_{n,i,j} &=& \text{the number of monotone $(2,n)$-trapezoids with }\\
& & \text{$1,2,\ldots, n$ in the last row and $i,j$ in the first row,}
\end{eqnarray*}
where a monotone $(2,n)$-trapezoid is a monotone triangle without its first row (see Fig.\ \ref{fig-illustr}), or again equivalently
\begin{eqnarray*}
A_{n,i,j} &=& \frac{1}{j-i+1}\times \Bigg(\text{the number of complete monotone triangles with}\\
& & \text{the numbers $i,j$ in the second row}\Bigg),
\end{eqnarray*}
since given the values $i,j$ in the second row, the number in the first row can be any integer in $[i,j]$ and has no influence on the numbers below the second row.

It is also possible to interpret the numbers $A_{n,i,j}$ directly in terms of the enumeration of classes of ASMs. Translating the above definitions to the language of ASMs, it is not difficult to see that for each triple $(i,j,k)$ such that $1 \le i\le k \le j \le n$ and $i<j$, we have that
\begin{eqnarray*}
%\frac{1}{j-i+1}
%(j-i+1)
A_{n,i,j} &=& \text{the number of ASMs $(m_{p,q})_{p,q=1}^n$ of order $n$ satisfying} \\
& & m_{1,k} = 1,\ \ m_{2,q}=0\text{ for $q\neq i,j,k$ and}\\&&  (m_{2,i}, m_{2,k}, m_{2,j}) = \left\{ \begin{array}{ll} (1,-1,1) & i<k<j, \\ (0,0,1) & i=k<j, \\ (1,0,0) & i<k=j, \end{array}\right.
\end{eqnarray*}
since the ASMs on the right-hand side correspond to monotone triangles with 
the number $k$ in the first row and the numbers $i,j$ in the second row.
Thus the numbers $A_{n,i,j}$ 
exactly enumerate ASMs based on their possible leading two rows; see Fig.~\ref{fig-asmtworows}.

\begin{figure}
\begin{center}
\begin{tabular}{ccc}
$ 
\begin{array}{l}
\begin{array}{cccccc}
\qquad \ i
\qquad \ \ k
\quad \ j
\end{array}\\ 
\left(\begin{array}{cccccc}
0 & 0 & 0 & 1 & 0 & 0 \\
0 & 1 & 0 & -1 & 1 & 0 \\
& & & \vdots & & 
\end{array}
\right)
\end{array}
$
&
$ 
\begin{array}{l}
\begin{array}{cccccc}
\qquad \qquad\ i\!\!=\!\!k
\  \ j
\end{array}\\ 
\left(\begin{array}{cccccc}
0 & 0 & 0 & 1 & 0 & 0 \\
0 & 0 & 0 & 0 & 1 & 0 \\
& & & \vdots & & 
\end{array}
\right)
\end{array}
$
&
$ 
\begin{array}{l}
\begin{array}{cccccc}
\qquad \ i
\quad\ \  k\!\!=\!\!j
\end{array}\\ 
\left(\begin{array}{cccccc}
0 & 0 & 0 & 1 & 0 & 0 \\
0 & 1 & 0 & 0 & 0 & 0 \\
& & & \vdots & & 
\end{array}
\right)
\end{array}
$
\\ (a) & (b) & (c)
\end{tabular}
\caption{The possible configurations for the first two rows of an ASM: (a) $i<k<j$. (b) $i=k<j$. (c) $i<k=j$. In all three cases the number of ASMs with the given first two rows is 
$A_{n,i,j}$.\label{fig-asmtworows}}
\end{center}
\end{figure}

Note that we have the special values
\begin{equation} \label{eq:specialvalues}
A_{n,i,n} = A_{n-1,i} \qquad (1\le i< n),
\end{equation}
and that the $A_{n,i,j}$'s satisfy the symmetry property
\begin{equation}\label{eq:symmetryeq}
A_{n,i,j} = A_{n,n+1-j,n+1-i} \qquad (1\le i<j\le n),
\end{equation}
which follows from the symmetry in the definition of ASMs with respect to horizontal reflection.

In what follows, we denote 
\begin{eqnarray*}
\binom{n}{k} &=& \left\{ \begin{array}{ll} \frac{n(n-1)(n-2)\ldots(n-k+1)}{k!} & k\ge 0, \\ 0 & \textrm{otherwise,}\end{array}\right. \\
\binom{n}{k}_+ &=& \left\{ \begin{array}{ll} \binom{n}{k} & n \ge 0, \\ 0 & \textrm{otherwise.}\end{array}\right.
\end{eqnarray*}
We can now formulate our results.

\begin{thm}[Linear equations]\label{mainresult}
Define
\begin{equation}
\ga_{i,j,p,q}=\begin{cases} (-1)^{i+q+1}\left({\binom {p-j+1} {q-i}}_+
-{\binom {p-j-1} {q-i-1}}_+\right)&\text{if }p \ge j,\\
0&\text{otherwise.}\end{cases}
\label{eq:gammas}
\end{equation}
For each $n\ge 1$, extend the triangular array of numbers $(A_{n,i,j})_{1\le i< j\le n}$ to a square matrix of numbers $(\hat{A}_{n,i,j})_{i,j=1}^n$ by setting
\begin{equation}\label{eq:extendednumbers}
\hat{A}_{n,i,j} = \left\{ \begin{array}{ll} A_{n,i,j} & i < j, \\
\sum_{1\le p<q\le n} \ga_{i,j,p,q} A_{n,p,q} & i \ge j.
\end{array}\right.
\end{equation}
Then the extended numbers $(\hat{A}_{n,i,j})_{i,j=1}^n$ satisfy the system of linear equations
\begin{equation} \label{eq:lineareqns}
\hat{A}_{n,i,j} = \sum_{p=i}^n \sum_{q=j}^n (-1)^{p+q} \binom{2n-i-2}{p-i}\binom{2n-j-2}{q-j} \hat{A}_{n,q,p}
\quad (1\le i,j \le n).
\end{equation}
\end{thm}

Figure \ref{fig-exampledoublyrefined} shows a table of the numbers $(\hat{A}_{n,i,j})_{i,j}$ for $n=5$. See Appendix A for a tabulation of these numbers for $3\le n\le 7$.
Note the following interesting special cases of \eqref{eq:extendednumbers}:
\begin{eqnarray}
\hat{A}_{n,n,j} &=& -\sum_{r=j}^{n-1} A_{n-1,r} \qquad (1\le j\le n), \label{eq:interestingid}\\
\hat{A}_{n,n,1} &=& -A_{n-1}, \nonumber \\
\hat{A}_{n,1,1} &=& \sum_{i=1}^{n-1} (-1)^{i+1} A_{n,i,i+1}, \nonumber \\
\hat{A}_{n,n,n} &=& 0. \nonumber
\end{eqnarray}

\begin{center}
\begin{figure}
$$
\left( \begin{array}{ccccc} 0&7&14&14&7\\-7&7&23&26&14\\
 -21&-2&16&23&14\\
 7&-7&-2&7&7\\
 -42&-35&-21&-7&0
\end{array}\right)
$$
\caption{The numbers $(\hat{A}_{5,i,j})_{i,j=1}^5$.
\label{fig-exampledoublyrefined}}
\end{figure}
\end{center}

We also prove the following result, which extends the symmetry property \eqref{eq:symmetryeq} in a surprising way.

\begin{thm}[Near-symmetry] \label{thm-almostsym}
The extended numbers $(\hat{A}_{n,i,j})_{i,j=1}^n$ satisfy the symmetry property
\begin{equation}\label{eq:symmetry}
\hat{A}_{n,i,j} = \hat{A}_{n,n+1-j,n+1-i}
\end{equation}
for all $1\le i,j\le n$, \dd{except} when $(i,j)=(n-1,1)$ or $(i,j)=(n,2)$, in which cases we have
\begin{equation}\label{eq:almostsym-exceptions}
\hat{A}_{n,n-1,1} = A_{n-2}, \qquad \hat{A}_{n,n,2} = A_{n-2}-A_{n-1}.
\end{equation}
\end{thm}

The near-symmetry equations \eqref{eq:symmetry} for $j\le i,\  (j,i)\notin \{(n-1,1),(n,2)\}$, are an additional set of linear relations satisfied by the $\hat{A}_{n,i,j}$'s, hence (by substituting the definition of $\hat{A}_{n,i,j}$ from \eqref{eq:extendednumbers} in both sides of \eqref{eq:symmetry}) by the original $A_{n,i,j}$'s. Some of these relations are quite simple -- for example, setting $i=j=1$ we get that
$ \hat{A}_{n,1,1} = \hat{A}_{n,n,n} = 0, $
which translates to
$$ \sum_{i=1}^{n-1} (-1)^{i+1} A_{n,i,i+1} = 0, $$
an identity which follows trivially from the symmetry of the $A_{n,i,j}$'s when $n$ is odd, but has interesting content when $n$ is even, and raises the question whether this equation has a direct combinatorial explanation. 

A natural question is whether the set of linear relations we found for the $A_{n,i,j}$'s, together with the ``boundary'' information \eqref{eq:specialvalues}, is sufficient to determine them. Strong numerical evidence suggests that the answer is positive.
\begin{conj}[Sufficiency] \label{conj-sufficiency}
The equations \eqref{eq:lineareqns}, together with the near-symmetry equations \eqref{eq:symmetry}, \eqref{eq:almostsym-exceptions}
and the special values \eqref{eq:specialvalues},
determine the numbers $\hat{A}_{n,i,j}$ uniquely.
\end{conj}

It is also natural to ask whether a closed-form expression analogous to the product formulas \eqref{eq:asmthm} and \eqref{eq:refinedasmthm} exists for the $A_{n,i,j}$'s. Assuming Conjecture \ref{conj-sufficiency}, because of Cramer's rule it should be possible to write a formula that expresses $A_{n,i,j}$ as a quotient of determinants. 
We also found the following formula, which is more explicit than a determinantal formula but still considerably more complicated than the formulas \eqref{eq:asmthm}, \eqref{eq:refinedasmthm}.

\begin{conj}[Explicit formula]\label{conj-explicitformula} 
For $(i,j)\notin\{ (n-1,1),(n,1),(n,2)\}$ we have:
\begin{eqnarray*}
\hat{A}_{n,i,j} &=& B_{n,i,j}(n+j-i-1 + P_{n,i,j} S_{n,i,j}),
 \quad \\ \ \\ & & \hspace{-70.0pt}\textrm{where we make the following definitions:} \\ \ \\ 
%A_n &=& \prod_{j=0}^{n-1} \frac{(3j+1)!}{(n+j)!}, \\
B_{n,i,j} &=& \frac{A_{n-1}}{(3n-5)!(n-2)!}\cdot \frac{(2n-2-i)!(2n-2-j)!(n+i-3)!(n+j-3)!}{(i-1)!(j-1)!(n-i)!(n-j)!}, \\
P_{n,i,j} &=& 2+2i+i^2-3j-i j+j^2-2n -2in + j n + n^2,  \\ %\ \\ \ \\
S_{n,i,j} &=& \sum_{k=\min(0, j-i)}^{\max(i-1,j-2)} \Big(X_{n,i,j}(k)-Y_{n,i,j}(k)\Big), \\ \ \\ %\ \\
X_{n,i,j}(k) &=& \left\{ \begin{array}{ll}
\frac{(-1)^{j+k+1}}{k-j+3-n}{3k-3j+4 \choose k} {2j+i-2k-5 \choose i-k-1} {i-2 \choose k-j+i}(i-1) \\
 \times
\left(3 H_{3j-2k-5}-3 H_{3j-3k-5} + 2H_{2j+i-2k-5} \right. \\  \ \left.- 2H_{2j-k-4}+H_{k-j+i}-H_{j-k-2}+\frac{1}{k-j+3-n} \right)\ \ 
& j-i\le k\le j-2,
\\ \ \\
\frac{ {3k-3j+4 \choose k} {2j+i-2k-5 \choose i-k-1} }{ {k-j+i \choose i-1} (k-j+3-n)}
& \textrm{otherwise,}
\end{array}\right.
\\ \  \\
Y_{n,i,j}(k) &=& \left\{ \begin{array}{ll}
\frac{(-1)^{i+k+1}}{k-j+3-n}{3k-3j+4 \choose k+i-j} {3j-2k-5 \choose j-k-1} {i-1 \choose k}(j-k-1) \\ \times
%\\ \qquad \times
\left(H_{3j-2k-5}-H_{2j-k-4}-H_k+H_{i-k-1} \right)
\ \ \ \quad\ \ %\ 
& 0\le k\le i-1,
\\ \ \\
\frac{ {3k-3j+4 \choose k+i-j } {3j-2k-5 \choose j-k-1 } (j-k-1) }{{k \choose i} (k-j+3-n) i}
& \textrm{otherwise,}
\end{array}\right.
\\
\ \\
%\ \\
H_m &=& \left\{ \begin{array}{ll} \sum_{d=1}^m \frac{1}{d} & m\ge 1, \\ \ \\ 0 & \textrm{otherwise.} \end{array}\right.
\end{eqnarray*}
For $(i,j)\in\{ (n-1,1),(n,1),(n,2)\}$ (where the values of $\hat{A}_{n,i,j}$ are known anyway) the formula above can be made to work by cancelling the zero of $P_{n,i,j}$ with the
appropriate singularity of $S_{n,i,j}$.
\end{conj}

We verified Conjectures \ref{conj-sufficiency} and \ref{conj-explicitformula} numerically up to $n=30$. 
See Appendix~B for details on how to compute the $\hat{A}_{n,i,j}$'s and verify the formulas using \texttt{RefinedASM}, the companion \texttt{Mathematica} software package to this paper. Note that $S_{n,i,j}$ defined above also has the following more succinct (but less explicit) representation:

$$ S_{n,i,j} = \lim_{j'\to j}  \Bigg[ \sum_{k=0}^{i-1} \left( \frac{\binom{3k-3j'+4}{k}\binom{2j'+i-2k-5}{i-1-k}}{\binom{k-j'+i}{i-1}(k-j'+3-n)}
-\frac{\binom{3k-3i+4}{k}\binom{2i+j'-2k-5}{i-1-k}(i-1-k)}{\binom{k-i+j'}{i}(k-i+3-n)i} \right)
\Bigg].
$$

\subsection{The $d$-refined enumeration}

The linear equations in Theorem \ref{mainresult} extend conjecturally to the $d$-refined enumeration of ASMs, for any $d\ge 1$. Once again, the formulation in terms of monotone triangles turns out to be the most convenient.

\begin{conj}[Linear equations for the $d$-refined case] \label{conj-lineareqs-drefined}
For $1\le d<n$ and $1\le i_1< i_2< \ldots < i_d\le n$, let $A_{n,i_1,i_2,\ldots,i_d}$ denote the number of
monotone triangles of order $n-d$ whose bottom row consists of the numbers in $\{1,2,\ldots,n\}\setminus \{i_1,i_2,\ldots,i_d\}$
arranged in increasing order
(if $d=n$ then we set $A_{n,1,2,\ldots,n}=1$). Then $(A_{n,i_1,\ldots,i_d})_{i_1<\ldots<i_d}$ can be extended into an array
$ (\hat{A}_{n,i_1,\ldots,i_d})_{1\le i_1,\ldots,i_d\le n} $ such that the equations
\begin{multline}
\label{eq:drefinedeqs}
\hat{A}_{n,i_1,\ldots,i_d} = (-1)^{nd}\sum_{j_1=i_1}^n \sum_{j_2=i_2}^n  \ldots \sum_{j_d=i_d}^n
(-1)^{j_1+\ldots+j_d} \\ \times\prod_{r=1}^d \binom{2n-i_r-d}{j_r-i_r} 
\hat{A}_{n,j_d,\ldots,j_2,j_1}
\qquad (1\le i_1,i_2,\ldots,i_d\le n)
\end{multline}
hold.
\label{conj-drefined}
\end{conj}

\subsection{A family of polynomial expansions}

We now formulate another result on the doubly-refined enumeration numbers (with a corresponding conjectural extension to the $d$-refined enumeration for $d\ge 3$) that is more conceptual than the linear equations \eqref{eq:lineareqns} described above, and whose proof will take up the main part of the paper. At the same time we also outline how this result will be used to derive the main results formulated above.

To introduce this result, we consider the function 
$$ \alpha_n(k_1, k_2, \ldots, k_n), $$
defined in \cite{fischer1}, that counts the number of monotone triangles with bottom row $(k_1,\ldots,k_n)$. This function is a polynomial with many useful symmetries. In \cite{fischer2} it was proved that the one-variable polynomial defined as a specialization of $\alpha_n$ by
$$ G_n(x) = \alpha_n(1,2,3,\ldots,n-1, n+x) $$
has the symmetry $G_n(x) = (-1)^{n-1} G_n(-2n-x)$ and that on the other hand, $G_n(x)$ is related to the refined enumeration numbers $A_{n,k}$ via the polynomial identity
\begin{equation}\label{eq:analogous}
G_n(x) = \sum_{k=1}^n A_{n,k} \binom{x+k-1}{k-1}.
\end{equation}
In other words, the $A_{n,k}$'s appear as the coefficients in the expansion of $G_n(x)$ in the linear basis $\left(\binom{x+k-1}{k-1}\right)_{k\ge 0}$ to the space of polynomials in $x$. By combining this with the symmetry property one gets a system of linear equations satisfied for each $n$ by the $A_{n,k}$'s, which was used in \cite{fischer2} to give a new proof of the refined alternating sign matrix theorem. 

We now extend this idea by considering the two-variable polynomial
\begin{equation}\label{eq:definitionG}
G_n(x,y) = \alpha_n(1,2,3,\ldots,n-2,n-1+x,n+y).
\end{equation}
As before, it will follow from the results of \cite{fischer2} that $G_n(x,y)$ satisfies a simple symmetry property. It turns out that $G_n(x,y)$ can, analogously to \eqref{eq:analogous}, be related to the (extended) doubly-refined enumeration numbers. We will prove the following identity.
\begin{thm}\label{thm-coefficients}
$$ G_n(x,y) = \sum_{i,j=1}^n \hat{A}_{n,i,j} \binom{x+i-1}{i-1}\binom{y+j}{j-1}. $$
\end{thm}
Again, this can be thought of as an expansion of the polynomial $G_n(x,y)$ in the linear basis
$ \{ \binom{x+i-1}{i-1}\binom{y+j}{j-1} \}_{i,j\ge 0} $ to the space of polynomials in $x$ and $y$.
The linear equations \eqref{eq:lineareqns} will follow relatively easily from a combination of this expansion with the symmetry property of $G_n(x,y)$. Along the way we gain enough information about the $\hat{A}_{n,i,j}$'s to also deduce Theorem \ref{thm-almostsym}.

Theorem \ref{thm-coefficients} and eq. \eqref{eq:analogous} also extend conjecturally to the $d$-refined enumerations for $d\ge3$, in the following form.

\begin{conj}
\label{conj-drefined}
For each $d\ge 1$, define the polynomials $G_n(x_1,\ldots,x_d)$ ($n\ge d$) by
$$ G_n(x_1,\ldots,x_d) = \alpha_n(1,2,\ldots,n-d,n-d+1+x_1,n-d+2+x_2,\ldots,n-1+x_{d-1},n+x_d). $$
Let $(F_{n,j_1,j_2,\ldots,j_d})_{j_1,\ldots,j_d=1}^n$ be the coefficients of  $G_n(x_1,\ldots,x_d)$
when expanded in the basis
$$ \left\{ 
\binom{x_1+j_1-1}{j_1-1}
\binom{x_2+j_2}{j_2-1}
\ldots
\binom{x_d+j_d+d-2}{j_d-1}
\right\}_{j_1,j_2,\ldots,j_d\ge 1}
$$
to the space of polynomials in $x_1,\ldots,x_d$, so that the equation
$$
G_n(x_1,\ldots,x_d) = \sum_{j_1,\ldots,j_d=1}^n F_{n,j_1,\ldots,j_d}
\prod_{r=1}^d \binom{x_r+j_r+r-2}{j_r-1}
$$
holds. Then we have, in the notation of Conjecture \ref{conj-lineareqs-drefined},
$$ F_{n,j_1,\ldots,j_d} = A_{n,j_1,\ldots, j_d} \qquad (1\le j_1<j_2<\ldots<j_d\le n). $$
\end{conj}

The coefficients $(F_{n,j_1,\ldots,j_d})_{j_1,\ldots,j_d}$ will give the ``extended'' enumeration numbers 
$(\hat{A}_{j_1,\ldots,j_d})_{j_1,\ldots,j_d}$ whose existence is the subject of Conjecture \ref{conj-lineareqs-drefined}.
In Section \ref{sectionproofmain} we will prove the following result, using the same ideas that were outlined above for how Theorem \ref{mainresult} can be deduced from Theorem {\ref{thm-coefficients}.

\begin{prop} \label{two-conjectures}
Conjecture \ref{conj-drefined} implies Conjecture \ref{conj-lineareqs-drefined}.
\end{prop}

\subsection{Acknowledgements}

We thank Christian Krattenthaler and Doron Zeilberger for helpful comments and suggestions.

\section{Preliminaries}

\subsection{The polynomial $\alpha_n(k_1,\ldots,k_n)$}

For integers $k_1< k_2<\ldots<k_n$, let $\alpha_n(k_1,\ldots,k_n)$ as before denote the number of monotone triangles with bottom row $(k_1,\ldots,k_n)$. From the definition of monotone triangles, the function $\alpha_n(k_1,\ldots,k_n)$ satisfies the recurrence relation
\begin{equation}\label{eq:alpharecurrence}
\alpha_n(k_1,\ldots,k_n) = \sum_{\small \begin{array}{c}j_1<\ldots<j_{n-1} \\k_1\le j_1\le k_2  \le \ldots \le j_{n-1}\le k_n\end{array}}\alpha_{n-1}(j_1,\ldots,j_{n-1}).
\end{equation}

It will be useful to note that if the $k_i$'s are only weakly increasing, i.e.,  $k_1\le k_2\le \ldots \le k_n$, then 
we can extend the definition of $\alpha_n(k_1,\ldots,k_n)$ to such a vector using the recurrence \eqref{eq:alpharecurrence}. This has the interpretation that this extended function counts the number of \dd{almost-monotone} triangles with a prescribed bottom row, where an almost-monotone triangle satisfies the same inequalities as a monotone triangle except that its bottom row is only required to be weakly  increasing.

In \cite[Th.~1]{fischer1} it was shown that one can extend the definition of $\alpha_n(k_1\ldots,k_n)$ even further to any vector of numbers, since it is a polynomial function in the variables $k_1,\ldots,k_n$ over $\mathbb{Q}$, which is of degree $n-1$ in every $k_i$. For convenience we will occasionally write simply $\alpha_n$ as a shorthand for the polynomial $\alpha_n(k_1\ldots,k_n)$.

The $A_{n,i,j}$'s can be expressed in terms of $\alpha_n$ as
$$
A_{n,i,j} = \alpha_{n-2}(1,2,\ldots,i-1,i+1,i+2\ldots,j-1,j+1,j+2,\ldots,n)
$$
for $1 \le i < j \le n$. 
By the above, the polynomial $G_n(x,y)$ defined in \eqref{eq:definitionG} can be expressed as a polynomial in $x$ and $y$ of degree no greater than $n-1$ in each of $x$ and $y$. Following the notation of Conjecture \ref{conj-drefined}, let $(F_{n,i,j})_{i,j=1}^n$ be
the coefficients of $G_n(x,y)$ in the expansion
\begin{equation}\label{eq:uniqueexp}
G_n(x,y) = \sum_{i=1}^{n} \sum_{j=1}^n F_{n,i,j} \binom{x+i-1}{i-1} \binom{y+j}{j-1},
\end{equation}
To reformulate Theorem \ref{thm-coefficients}, we need to prove that $F_{n,i,j}=\hat{A}_{n,i,j}$ for all $n\ge 3$ and $1\le i,j\le n$.

\subsection{Two lemmas}

If $x$ is a variable, let $E_x$ denote the right-shift operator in the variable $x$, acting on polynomials in $x$ by
$$ (E_x h)(x) = h(x+1), $$
and let $\Delta_x$ denote the (right-)differencing operator in $x$, defined by
$$ (\Delta_x h)(x) = h(x+1)-h(x) = ((E_x-\textrm{Id})h)(x). $$
We need the following lemma, which appeared in equivalent form in \cite[Lemma~1]{fischer2}. 
\begin{lemma} 
\label{sym}
Let $P(X_1,\ldots,X_n)$ be a symmetric polynomial in $X_1,\ldots, X_n$. Then
$$
P(\Delta_{k_1},\ldots,\Delta_{k_n}) \alpha_n(k_1,\ldots,k_n) = P(0,\ldots,0) \alpha_n(k_1,\ldots,k_n).
$$
\end{lemma}

Lemma \ref{sym} is used to show the next lemma, which will be applied twice in the following. It involves the $p$--th 
elementary symmetric function, denoted by 
$$
e_p(X_1,\ldots,X_n): = \sum_{1 \le i_1 < i_2 \ldots < i_p \le n} X_{i_1} X_{i_2} \dots X_{i_p}
$$ 
(for $p=0$, we take $e_0(X_1,\ldots,X_n)=1$).
As before, the notation $a_1,\ldots, \widehat{a_i},\ldots, a_m$ for a list of objects is used to denote the list with the element $a_i$ omitted from it.
\begin{lemma} 
\label{twice}
Let $z$ be a non-negative integer, $m \ge 1$ and $1 \le r \le m$. Then 
$$
E^z_{k_r} \alpha_m(k_1,\ldots,k_m) 
= (-1)^z \sum_{p=0}^{z} 
\binom{-m}{z-p} e_p(E_{k_1},\ldots,\widehat{E_{k_r}},\ldots,E_{k_m}) 
\alpha_m(k_1,\ldots,k_m).
$$
\end{lemma}

\begin{proof}
For any $q>0$ we have $e_q(0,\ldots,0)=0$ and
$$
e_q(\Delta_{k_1},\ldots,\widehat{\Delta_{k_r}},\ldots,\Delta_{k_m}) =
e_q(\Delta_{k_1},\ldots,\Delta_{k_m}) - 
\Delta_{k_r}e_{q-1}(\Delta_{k_1},\ldots,\widehat{\Delta_{k_r}},\ldots,\Delta_{k_m}).
$$
It follows by induction, using Lemma \ref{sym}, that
$$
\Delta_{k_r}^q \alpha_m(k_1,\ldots,k_m) = (-1)^q 
e_q(\Delta_{k_1},\ldots,\widehat{\Delta_{k_r}},\ldots,\Delta_{k_m})
\alpha_m(k_1,\ldots,k_m).
$$
Therefore
\begin{eqnarray*}
E^z_{k_r} \alpha_m(k_1,\ldots,k_m)
&=& 
(\Delta_{k_r} + \id)^z 
\alpha_m(k_1,\ldots,k_m)
\qquad\qquad\qquad\qquad\qquad\qquad
\end{eqnarray*}

\vspace{-20.0pt}
\begin{eqnarray}
\qquad &=& \sum_{q=0}^{z} \binom{z}{q} \Delta_{k_r}^q  \alpha_m(k_1,\ldots,k_m) \nonumber \\
&=& \sum_{q=0}^z \binom{z}{q} (-1)^q e_q(\Delta_{k_1},\ldots,\widehat{\Delta_{k_r}},\ldots,\Delta_{k_m})
\alpha_m(k_1,\ldots,k_m). 
\label{ez}
\end{eqnarray}
Since 
\begin{eqnarray*}
e_q(\Delta_{k_1},\ldots,\widehat{\Delta_{k_r}},\ldots,\Delta_{k_m}) &=& 
e_q(E_{k_1}-\id,\ldots,\widehat{E_{k_r}-\id},\ldots,E_{k_m}-\id) \\ &
=& \sum_{p=0}^{q} \binom{m-1-p}{q-p} (-1)^{q-p} e_p(E_{k_1},\ldots,\widehat{E_{k_r}},\ldots,E_{k_m}),
\end{eqnarray*}
the right-hand side of \eqref{ez} is equal to 
\begin{multline*}
\sum_{q=0}^z \binom{z}{q} (-1)^q \sum_{p=0}^{q} \binom{m-1-p}{q-p} (-1)^{q+p} 
e_p(E_{k_1},\ldots,\widehat{E_{k_r}},\ldots,E_{k_m}) \alpha_m(k_1,\ldots,k_m) \\
= \sum_{p=0}^{z} (-1)^p \sum_{q=p}^{z} \binom{z}{q} \binom{m-1-p}{q-p} 
e_p(E_{k_1},\ldots,\widehat{E_{k_r}},\ldots,E_{k_m}) 
\alpha_{m}(k_1,\ldots,k_m).
\end{multline*}
Using the Chu-Vandermonde summation formula 
$$ \binom{a+b}{c} = \sum_k \binom{a}{k} \binom{b}{c-k} $$
we now get
\begin{eqnarray*}
(-1)^p \sum_{q=p}^{z} \binom{z}{q} \binom{m-1-p}{q-p} 
&=& (-1)^p \sum_{q=p}^{z} \binom{z}{z-q} \binom{m-1-p}{q-p} \\
&=& (-1)^p \sum_{r=0}^{z-p} \binom{z}{z-p-r} \binom{m-1-p}{r} 
\\&=& (-1)^p \binom{z+m-1-p}{z-p} = (-1)^z \binom{-m}{z-p}
\end{eqnarray*}
and the assertion follows.
\end{proof}

\section{A first representation for $F_{n,i,j}$\label{fz}}

For $p\ge 0$ and $1\le i\le n$ we set 
$$
Z(n,p,i): =  \Big(e_{p}(E_{k_1},E_{k_2},\ldots,E_{k_{n-2}}) \alpha_{n-1}\Big)(1,2,\ldots,i-1,i+1,i+2,\ldots,n),
$$
and also define $Z(n,p,0)=0$.
The proof of Theorem \ref{thm-coefficients} will proceed in several stages. The first of these will be to deduce
a formula for the coefficients $F_{n,i,j}$ in terms of the numbers $\{Z(n,p,i)\}_p$.

In the following, we assume $n \ge 3$. We apply Lemma~\ref{twice} with the parameters $m=n-1$, $r=n-1$, $z=y$ and 
$(k_1,\ldots,k_{n-1})=(1,2,\ldots,i-1,i+1,i+2,\ldots,n)$, to get that 
\begin{eqnarray}
C(n,i,y) &:=& \alpha_{n-1}(1,2,\ldots,i-1,i+1,i+2,\ldots,n-1,n+y)  \nonumber
\\ &=&(-1)^y \sum_{p=0}^{y} Z(n,p,i) \binom{-n+1}{y-p}
\label{2}
\end{eqnarray}
for all non-negative integers $y$.

We now consider
$
\Big( e_{q}(E_{k_1},E_{k_2},\ldots,E_{k_{n-2}}) \alpha_{n}\Big)(1,2,\ldots,n-1,n+y).
$
If $1 \le i_1 < i_2 < \dots < i_q \le n-2$ then 
\begin{multline}
\left( E_{k_{i_1}} E_{k_{i_2}} \dots E_{k_{i_q}}  \alpha_{n}\right)(1,2,\ldots,n-1,n+y) \\
= \left(E_{k_{i_1}} \alpha_n\right)(1,2,\ldots,n-1,n+y) \qquad\qquad\qquad\qquad\qquad\qquad\quad\ \ \\
= \alpha_n(1,2,\ldots,i_1-1,i_1+1,i_1+1,i_1+2,\ldots,n-1,n+y).\qquad\qquad
\label{eq:followsfrom}
\end{multline}
This follows from the recurrence equation 
\eqref{eq:alpharecurrence}, which implies that the quantity on the right-hand side does not change if some of the numbers in position $i_1+1,i_1+2,\ldots,n-2$ in the argument of $\alpha_n$ are increased by one. 
Moreover, again by the recurrence \eqref{eq:alpharecurrence}, and using the identity
$\binom{n}{k}=(-1)^k \binom{k-n+1}{k}$, we get
\begin{samepage}
\begin{equation*}
\alpha_n(1,2,\ldots,i_1-1,i_1+1,i_1+1,i_1+2,\ldots,n-1,n+y) \qquad\qquad\qquad\qquad
\end{equation*}

\vspace{-24.0pt}
\begin{eqnarray}
\qquad\ \  & =& \sum_{i=1}^{i_1} \sum_{y_1=0}^{y} C(n,i,y_1)  \nonumber \\
&=& \sum_{i=1}^{i_1} \sum_{y_1=0}^{y} \sum_{p=0}^{y_1} (-1)^{y_1} Z(n,p,i) \binom{-n+1}{y_1-p}  \nonumber \\
&=& \sum_{i=1}^{i_1} \sum_{p=0}^{y} \sum_{y_1=p}^{y} (-1)^{p} Z(n,p,i) \binom{y_1-p+n-2}{n-2} \nonumber \\
&=& \sum_{i=1}^{i_1} \sum_{p=0}^{y} (-1)^{p} Z(n,p,i) \sum_{y_1=p}^{y} \left( \binom{y_1-p+n-1}{n-1}
- \binom{y_1-p+n-2}{n-1} \right)  \nonumber \\
&=& \sum_{i=1}^{i_1} \sum_{p=0}^y (-1)^{p} Z(n,p,i) \binom{y-p+n-1}{n-1}.
\label{eff}
\end{eqnarray}
\end{samepage}

\noindent
Now divide into two cases, $q=0$ and $q\neq 0$. If $q\neq 0$,
then, since there are $\binom{n-2-i_1}{q-1}$ subsets of $\{1,2,\ldots,n-2\}$ with $q$ elements and 
minimal element $i_1$, it follows from \eqref{eq:followsfrom} and \eqref{eff} that
\begin{samepage}
\begin{equation*}
 \Big( e_{q}(E_{k_1},E_{k_2},\ldots,E_{k_{n-2}}) \alpha_{n}\Big)(1,2,\ldots,n-1,n+y)
 \qquad\qquad\qquad\qquad\qquad\qquad\qquad
\end{equation*}

\vspace{-20.0pt}
\begin{eqnarray}
\ \ &=& \sum_{i_1=1}^{n-2} \binom{n-2-i_1}{q-1} \sum_{i=1}^{i_1} \sum_{p=0}^y (-1)^{p} Z(n,p,i) \binom{y-p+n-1}{n-1}  \nonumber \\
&=& \sum_{i=1}^{n-2} \sum_{i_1=i}^{n-2} \binom{n-2-i_1}{q-1} \sum_{p=0}^y (-1)^{p} Z(n,p,i) \binom{y-p+n-1}{n-1}  \nonumber \\
&=& \sum_{i=1}^{n-2} \sum_{i_1=i}^{n-2} \left( \binom{n-1-i_1}{q} - \binom{n-2-i_1}{q} \right) 
\sum_{p=0}^y (-1)^{p} Z(n,p,i) \binom{y-p+n-1}{n-1}  \nonumber \\
&=& \sum_{i=1}^{n-1} \sum_{p=0}^y (-1)^{p} Z(n,p,i) \binom{y-p+n-1}{n-1} \binom{n-1-i}{q}.
\label{eq:unifying1}
\end{eqnarray}
\end{samepage}
If $q=0$, then using \eqref{eq:alpharecurrence} and using the equality between the second expression and the last expression
in \eqref{eff}, we get
\begin{multline}
\alpha_{n}(1,2,\ldots,n-1,n+y) = \sum_{i=1}^{n-1} \sum_{p=0}^{y} C(n,i,p) +
\alpha_{n-1}(1,2,\ldots,n-1) \\
 =  \sum_{i=1}^{n-1} \sum_{p=0}^y (-1)^{p} Z(n,p,i) \binom{y-p+n-1}{n-1} + A_{n-1}. \label{eq:unifying2}
\end{multline}
So, unifying \eqref{eq:unifying1} and \eqref{eq:unifying2} we can write
\begin{multline}
 \Big( e_{q}(E_{k_1},E_{k_2},\ldots,E_{k_{n-2}}) \alpha_{n}\Big)(1,2,\ldots,n-1,n+y)
 \\=
\sum_{i=1}^{n-1} \sum_{p=0}^y (-1)^{p} Z(n,p,i) \binom{y-p+n-1}{n-1} \binom{n-1-i}{q}
+ [q=0] A_{n-1},
\label{form}
\end{multline}
where $[\text{statement}]=1$ if statement is true and $[\text{statement}]=0$ otherwise.
As a consequence of \eqref{form} we get that
\begin{multline}
\label{3}
 \Big( e_{q}(E_{k_1},E_{k_2},\ldots,E_{k_{n-2}},E_{k_n}) \alpha_{n}\Big)(1,2,\ldots,n-1,n+y)  \\
= \Big( \Big(e_{q}(E_{k_1},\ldots,E_{k_{n-2}}) + e_{q-1}(E_{k_1},\ldots,E_{k_{n-2}}) E_{k_n}\Big)
\alpha_{n}\Big)(1,2,\ldots,n-1,n+y)  \\
= \sum_{i=1}^{n-1} \sum_{p=0}^y (-1)^{p} Z(n,p,i) \binom{y-p+n-1}{n-1} \binom{n-1-i}{q}  
\qquad\qquad\qquad\qquad\qquad\ \ 
\\
+ \sum_{i=1}^{n-1} \sum_{p=0}^{y+1} (-1)^{p} Z(n,p,i) \binom{y-p+n}{n-1} \binom{n-1-i}{q-1}
+ [q \in \{0,1\}] A_{n-1}. \quad
\end{multline}

Now use Lemma~\ref{twice} again, with $m=n$, $r=n-1$, $z=x$ and $(k_1,\ldots,k_n)=(1,2,\ldots,n-1,n+y)$. This gives the following representation for the polynomial $G_n(x,y)$ defined in \eqref{eq:definitionG}:
\begin{multline*}
G_n(x,y)= \alpha_{n}(1,2,\ldots,n-2, n-1+x,n+y) \\
= (-1)^x \sum_{q=0}^{x} \binom{-n}{x-q} \Big( e_{q}(E_{k_1},E_{k_2},\ldots,E_{k_{n-2}},E_{k_n}) 
\alpha_{n}\Big)(1,2,\ldots,n-1,n+y).
\end{multline*}
Combining this with \eqref{3}, and using Chu-Vandermonde summation again, we get that
\begin{samepage}
\begin{equation*}
 G_n(x,y)
 \qquad\qquad\qquad\qquad\qquad\qquad\qquad\qquad\qquad\qquad
 \qquad\qquad\qquad\qquad \qquad\qquad
\end{equation*}

\vspace{-26.0pt}
\begin{eqnarray*}
 &=& (-1)^x \sum_{q=0}^{x} [q \in \{0,1\}] \binom{-n}{x-q} A_{n-1} 
+ (-1)^x \sum_{q=0}^{x} \sum_{i=1}^{n-1} \sum_{p=0}^{y+1} (-1)^{p} Z(n,p,i)
 \binom{-n}{x-q} \\
& & \qquad \times \left( \binom{y-p+n-1}{n-1} \binom{n-1-i}{q} +  \binom{y-p+n}{n-1} \binom{n-1-i}{q-1} \right) 
\\ &=& \binom{x+n-2}{n-2} A_{n-1} 
+ \sum_{i=1}^{n-1} \sum_{p=0}^{y+1} (-1)^{p} Z(n,p,i) \\
& & \qquad \times
 \left( \binom{y-p+n-1}{n-1} \binom{x+i}{i} -  \binom{y-p+n}{n-1} \binom{x+i-1}{i} \right).
\end{eqnarray*}
\end{samepage}
\noindent This identity is valid for all integer $x,y\ge 0$. Therefore the two polynomials in $x$ and $y$ are in fact identical.

We want to expand this polynomial in terms of $\binom{x+i-1}{i-1} \binom{y+j}{j-1}$. Observe that 
\begin{samepage}
\begin{equation*}
\binom{y-p+n-1}{n-1}\binom{x+i}{i} - \binom{y-p+n}{n-1}\binom{x+i-1}{i}
\qquad\qquad\qquad\qquad\qquad\ \ 
\end{equation*}

\vspace{-16.0pt}
\begin{eqnarray*}
&=&
\Bigg( \binom{y-p+n}{n-1}-\binom{y-p+n-1}{n-2}\Bigg) \binom{x+i}{i}
 - \binom{y-p+n}{n-1}\binom{x+i-1}{i} \\&=&
 \binom{y-p+n}{n-1}\Bigg(\binom{x+i}{i}-\binom{x+i-1}{i}\Bigg)
 -\binom{y-p+n-1}{n-2}\binom{x+i}{i}
 \\ &=& \binom{y-p+n}{n-1} \binom{x+i-1}{i-1} -  \binom{y-p+n-1}{n-2} \binom{x+i}{i},
\end{eqnarray*}
\end{samepage}

\noindent
and therefore
\begin{eqnarray*}
G_n(x,y)
 &=& \binom{x+n-2}{n-2} A_{n-1}  
 + \sum_{i=1}^{n-1} \sum_{p=0}^{n-2} (-1)^{p} Z(n,p,i) 
 \\ & & \qquad%\qquad\qquad\qquad 
 \times \left( \binom{y-p+n}{n-1} \binom{x+i-1}{i-1} -  \binom{y-p+n-1}{n-2} \binom{x+i}{i} \right).
\end{eqnarray*}
By Chu-Vandermonde summation we have 
$$
\binom{y-p+m}{m-1} = \sum_{j=1}^{m} (-1)^{j+m} \binom{p}{m-j} \binom{y+j}{j-1}
$$
and therefore
\begin{multline*}
 G_n(x,y)
 =   \binom{x+n-2}{n-2} A_{n-1}  \\
 + \sum_{i=1}^{n-1} \sum_{p=0}^{n-2} (-1)^{p} Z(n,p,i) \binom{x+i-1}{i-1} \sum_{j=1}^{n} (-1)^{j+n} \binom{p}{n-j} 
 \binom{y+j}{j-1} \ \ \ \ \ \\
 + \sum_{i=1}^{n} \sum_{p=0}^{n-2} (-1)^{p} Z(n,p,i-1) \binom{x+i-1}{i-1} \sum_{j=1}^{n} (-1)^{j+n} \binom{p}{n-1-j} \binom{y+j}{j-1}, 
 \end{multline*}
(recall that $Z(n,p,0)=0$).
The coefficient of $\binom{x+i-1}{i-1} \binom{y+j}{j-1}$ in this expression is
\begin{eqnarray}
F_{n,i,j} &=& \sum_{p=0}^{n-2} (-1)^{p+j+n} \left( [i \not= n] Z(n,p,i) \binom{p}{n-j} + Z(n,p,i-1) \binom{p}{n-1-j} \right) \nonumber \\
& & \qquad + [i=n-1][j=1] A_{n-1}. \label{eq:coefficientis}
\end{eqnarray}

\section{A second representation for $F_{n,i,j}$\label{za}}

Now define
\begin{equation*}
W(n,i,j) := \sum_{p=0}^{n-2} (-1)^{p+j+n} Z(n,p,i) \binom{p}{n-j}.
\end{equation*}
By \eqref{eq:coefficientis}, $F_{n,i,j}$ can be expressed in terms of the $W(n,s,t)$'s as
\begin{equation}
F_{n,i,j} = [i\neq n]W(n,i,j)-W(n,i-1,j+1)+[i=n-1][j=1]A_{n-1}. \label{eq:Fexpressedas}
\end{equation}

We derive a formula for $W(n,i,j)$. For $1 \le i' \le i \le n-1$, set 
\begin{multline*}
 Q(n,p,i',i)= \\ \Big(
e_p(E_{k_1},\ldots,E_{k_{i'-2}},E_{k_i},\ldots,E_{k_{n-2}}) \alpha_{n-1}\Big)(1,2\ldots,i'-1,i'+1,i'+2,\ldots,n),
\end{multline*}
and observe that 
\begin{samepage}
\begin{equation*}
Z(n,p,i) = \Big( e_{p}(E_{k_1},E_{k_2},\ldots,E_{k_{n-2}}) \alpha_{n-1}\Big)(1,2,\ldots,\widehat{i},\ldots,n)
\qquad\qquad\qquad\qquad
\end{equation*}

\vspace{-22.0pt}
\begin{eqnarray*}
&=& \Big( e_{p}(E_{k_1},\ldots,\widehat{E_{k_{i-1}}},\ldots,E_{k_{n-2}}) \alpha_{n-1}\Big)(1,2,\ldots,\widehat{i},\ldots,n)\\
& & \qquad + \Big( e_{p-1}(E_{k_1},\ldots,\widehat{E_{k_{i-1}}},\ldots,E_{k_{n-2}}) \alpha_{n-1}\Big)(1,2,\ldots,\widehat{i-1},\ldots,n) \\
& =& Q(n,p,i,i) \\
& & \qquad + \Big( e_{p-1}(E_{k_1},\ldots,\widehat{E_{k_{i-2}}},\widehat{E_{k_{i-1}}},\ldots,E_{k_{n-2}}) \alpha_{n-1}
\Big)(1,2,\ldots,\widehat{i-1}\ldots,n) \\
& & \qquad
+ \Big( e_{p-2}(E_{k_1},\ldots,\widehat{E_{k_{i-2}}},\widehat{E_{k_{i-1}}},\ldots,E_{k_{n-2}}) \alpha_{n-1}\Big)(1,2,\ldots,\widehat{i-2}\ldots,n) \\
&=& \ldots = \sum_{i'=1}^{i} Q(n,p-i+i',i',i)
\end{eqnarray*}
\end{samepage}
(for $i=0$, $Z(n,p,i)$ is 0 and therefore also trivially equal to the expression on the right-hand side).
Therefore,
\begin{equation}
\label{toshow}
W(n,i,j)=\sum_{p=0}^{n-2} \sum_{i'=1}^{i} (-1)^{p+j+n} \binom{p}{n-j} Q(n,p-i+i',i',i).
\end{equation}

Next we will see how $Q(n,p,i',i)$ can be expressed in terms of the $A_{n,a,b}$'s. For each 
$P \subseteq \{1,2,\ldots,i'-2\} \cup \{i,i+1,\ldots,n-2\}$, set
$$Q(n,P,i')  = \Big( \Big( \prod_{r \in P} E_{k_r} \Big) \alpha_{n-1}\Big)(1,2,\ldots,i'-1,i'+1,i'+2,\ldots,n). $$
Then $Q(n,p,i',i)$ is the sum of $Q(n,P,i')$ over all such sets $P$ of size $p$. We evaluate $Q(n,P,i')$, separating into 4 cases
according as whether $P_1 := P\cap \{1,\ldots,i'-2\}$ is nonempty and whether $P_2 := P\cap\{i,i+1,\ldots,n-2\}$ is nonempty. In each of these cases, we let $s_j$ denote the minimal element of $P_j$ ($j=1,2$) if $P_j$ is nonempty.

\paragraph{Case 1: $P_1$ and $P_2$ are both nonempty.} In this case, from the recurrence \eqref{eq:alpharecurrence}
it follows using an argument similar to the one used previously that
\begin{eqnarray*}
\nonumber
Q(n,P,i')&=&
\alpha_{n-1}(1,2,\ldots,s_1-1,s_1+1,s_1+1,s_1+2,\ldots,\\ & & \qquad\ \ \ \ 
i'-1,i'+1,i'+2,\ldots,s_2,s_2+2,s_2+2,s_2+3\ldots,n).
\nonumber
\\ &=&\sum_{1 \le a \le s_1,\ i' \le b \le s_2+1} A_{n,a,b}.
\label{case1}
\end{eqnarray*}
 
\paragraph{Case 2: $P_1$ is nonempty, $P_2$ is empty.} In this case, by a similar argument we have that
\begin{eqnarray*}
\nonumber
Q(n,P,i') &=&\alpha_{n-1}(1,2,\ldots,s_1-1,s_1+1,s_1+1,s_1+2,\ldots, \\ & & \qquad\qquad\qquad\qquad\qquad\qquad 
i'-1,i'+1,i'+2,\ldots,n) \nonumber \\
&=& \sum_{1 \le a \le s_1,\ i' \le b \le n} A_{n,a,b}.
\label{case2}
\end{eqnarray*}

\paragraph{Case 3: $P_1$ is empty, $P_2$ is nonempty.} We have
\begin{eqnarray*}
\nonumber
Q(n,P,i') &=& 
\alpha_{n-1}(1,2,\ldots,i'-1,i'+1,i'+2,\ldots,\\ & & \qquad\qquad\qquad
s_2,s_2+2,s_2+2,s_2+3,\ldots,n)
\nonumber \\ &=& \sum_{1 \le a \le i' \le b \le s_2+1,\ a < b} A_{n,a,b}.
\label{case3}
\end{eqnarray*}

\paragraph{Case 4: $P=P_1=P_2=\emptyset$.} In this case, $p=|P|=0$ and
\begin{eqnarray*}
Q(n,P,i')&=& \alpha_{n-1}(1,2,\ldots,\widehat{i'},\ldots,n) = A_{n,i'} \\ &=&
\sum_{1 \le a \le i' \le b \le n,\ a < b} A_{n,a,b}.
\end{eqnarray*}

Note that the number of sets $P$ of size $p$ with given values of $s_1, s_2$ for which Case 1 holds is $\binom{i'+n-4-s_1-s_2}{p-2}$; the number of sets of size $p$ with given $s_1$ for which Case 2 holds is
$\binom{i'-2-s_1}{p-1}$; Case 3 holds for $\binom{n-2-s_2}{p-1}$ sets of size $p$ with given $s_2$, and Case 4 holds for one set if $p=0$, or for no sets otherwise. Putting all the above information together, we get the following formula for $Q(n,p,i',i)$.
\begin{eqnarray*}
Q(n,p,i',i) &=& 
\sum_{s_1=1}^{i'-2} \sum_{s_2=i}^{n-2} \binom{i'+n-4-s_1-s_2}{p-2} \sum_{1 \le a \le s_1,\ i' \le b \le s_2+1} 
A_{n,a,b} \\
& &  + \sum_{s_1=1}^{i'-2} \binom{i'-2-s_1}{p-1} \sum_{1 \le a \le s_1,\ i' \le b \le n} A_{n,a,b} 
\\ & &   + \sum_{s_2=i}^{n-2} \binom{n-2-s_2}{p-1} \sum_{1 \le a \le i' \le b \le s_2+1,\ a < b} A_{n,a,b} \\
& &  + [p=0] \sum_{1 \le a \le i' \le b \le n,\ a < b} A_{n,a,b}.
\end{eqnarray*}
We change the order  of summation, split up some terms into two parts, and adopt the convention that $A_{n,i,i}=0$,
to see that
\begin{eqnarray*}
Q(n,p,i',i)&\!\!\!=\!\!\!&
\sum_{a=1}^{i'-2} \sum_{b=i'}^{i} \sum_{s_1=a}^{i'-2}
\sum_{s_2=i}^{n-2} \binom{i'+n-4-s_1-s_2}{p-2} A_{n,a,b} 
\\ & &+ \sum_{a=1}^{i'-2} \sum_{b=i+1}^{n-1} \sum_{s_1=a}^{i'-2} \sum_{s_2=b-1}^{n-2}
\binom{i'+n-4-s_1-s_2}{p-2} A_{n,a,b} \\
& & + \sum_{a=1}^{i'-2} \sum_{b=i'}^{n} \sum_{s_1=a}^{i'-2} \binom{i'-2-s_1}{p-1}
A_{n,a,b} 
+ \sum_{a=1}^{i'} \sum_{b=i'}^{i} \sum_{s_2=i}^{n-2} \binom{n-2-s_2}{p-1}
A_{n,a,b} \\ & &
+ \sum_{a=1}^{i'} \sum_{b=i+1}^{n-1} \sum_{s_2=b-1}^{n-2} \binom{n-2-s_2}{p-1}
A_{n,a,b} + 
[p=0] \sum_{a=1}^{i'} \sum_{b=i'}^{n} A_{n,a,b}.
\end{eqnarray*}
Applying the telescopic summation 
\begin{equation}
\label{sum}
\sum_{z=x}^{y} \binom{c-z}{m} = \sum_{z=r}^{s} \left( \binom{c-z+1}{m+1} - \binom{c-z}{m+1} \right) 
= \binom{c-r+1}{m+1} - \binom{c-s}{m+1}
\end{equation}
eliminates one summation operator from each term, and after some minor rearrangement (using the fact that $\binom{0}{p}=[p=0]$) brings this to the form
\begin{eqnarray*}
Q(n,p,i',i) &\!\!\!=\!\!\!& \sum_{a=1}^{i'-2} \sum_{b=i'}^{i} \sum_{s_2=i}^{n-2} \left(
  \binom{i'+n-3-a-s_2}{p-1} - \binom{n-2-s_2}{p-1} \right) A_{n,a,b} \\ & &
+\sum_{a=1}^{i'-2} \sum_{b=i+1}^{n-1} \sum_{s_2=b-1}^{n-2} \left(
  \binom{i'+n-3-a-s_2}{p-1} - \binom{n-2-s_2}{p-1} \right) A_{n,a,b} \\ & &
+\sum_{a=1}^{i'-2} \sum_{b=i'}^{n} \left(  \binom{i'-1-a}{p} - \binom{0}{p}
\right) A_{n,a,b} \\
& & + \sum_{a=1}^{i'} \sum_{b=i'}^{i} \binom{n-1-i}{p} A_{n,a,b}
+ \sum_{a=1}^{i'} \sum_{b=i+1}^{n-1} \binom{n-b}{p} A_{n,a,b}
\\ & & +[p=0] \sum_{a=1}^{i'} A_{n,a,n}.
\end{eqnarray*}
Applying \eqref{sum} once again to the first two terms then gives
\begin{samepage}
$$ Q(n,p,i',i) 
\qquad\qquad\qquad\qquad\qquad\qquad\qquad\qquad\qquad\qquad\qquad\qquad
\qquad\qquad\qquad
$$

\vspace{-38.0pt}
\begin{eqnarray*}
\\ &=& \sum_{a=1}^{i'-2} \sum_{b=i'}^{i} \left( \binom{i'+n-2-a-i}{p} - \binom{i'-a-1}{p} \right) A_{n,a,b} 
\\ & &+ \sum_{a=1}^{i'-2} \sum_{b=i'}^{i} \left( - \binom{n-1-i}{p} + \binom{0}{p} \right) A_{n,a,b} \\
& & + \sum_{a=1}^{i'-2} \sum_{b=i+1}^{n-1} \left( \binom{i'+n-1-a-b}{p} - \binom{i'-a-1}{p} \right) A_{n,a,b} 
\\ & &+ \sum_{a=1}^{i'-2} \sum_{b=i+1}^{n-1} \left( - \binom{n-b}{p} + \binom{0}{p} \right) A_{n,a,b} 
\\ & & + \sum_{a=1}^{i'-2} \sum_{b=i'}^n \left( \binom{i'-1-a}{p} - \binom{0}{p} \right) A_{n,a,b} \\ & &
+ \sum_{a=1}^{i'} \sum_{b=i'}^{i} \binom{n-1-i}{p} A_{n,a,b} + \sum_{a=1}^{i'} \sum_{b=i+1}^{n-1} \binom{n-b}{p} A_{n,a,b} 
 + [p=0] \sum_{a=1}^{i'} A_{n,a,n}, 
\end{eqnarray*}
\end{samepage}

\begin{samepage}
\noindent which furthermore simplifies to
$$
Q(n,p,i',i)
\qquad\qquad\qquad\qquad\qquad\qquad\qquad\qquad\qquad\qquad\qquad\qquad
\qquad\qquad\qquad
$$

\vspace{-24.0pt}
\begin{eqnarray*}
&=& \sum_{a=1}^{i'-1} \sum_{b=i'}^{i} \binom{i'+n-2-a-i}{p} A_{n,a,b} 
+ \sum_{a=1}^{i'-1} \sum_{b=i+1}^{n} \binom{i'+n-1-a-b}{p} A_{n,a,b} \\
& & + \sum_{b=i'+1}^{i} \binom{n-1-i}{p} A_{n,i',b} 
+ \sum_{b=i+1}^{n} \binom{n-b}{p} A_{n,i',b}.
\end{eqnarray*}
\end{samepage}

\noindent We now substitute this formula for $Q(n,p,i',i)$ in \eqref{toshow} and interchange summation operators to get
\begin{eqnarray*}
W(n,i,j) &=& 
\sum_{p=0}^{n-2} \sum_{a=1}^{i-1} \sum_{b=a+1}^{i} \sum_{i'=a+1}^{b} (-1)^{p+j+n} \binom{p}{n-j} 
\binom{i'+n-2-a-i}{n-a-p-2} A_{n,a,b}  \\
& & + 
\sum_{p=0}^{n-2} \sum_{a=1}^{i-1} \sum_{b=i+1}^{n} \sum_{i'=a+1}^{i} (-1)^{p+j+n} \binom{p}{n-j} 
\binom{i'+n-1-a-b}{n-a-b-p+i-1} A_{n,a,b} \\
& & + \sum_{p=0}^{n-2} \sum_{b=2}^i \sum_{i'=1}^{b-1} (-1)^{p+j+n} \binom{p}{n-j} \binom{n-1-i}{p-i+i'} A_{n,i',b} \\
& &+ \sum_{p=0}^{n-2} \sum_{b=i+1}^{n} \sum_{i'=1}^{i} (-1)^{p+j+n} \binom{p}{n-j} \binom{n-b}{p-i+i'} A_{n,i',b}.
\end{eqnarray*}
Then, applying \eqref{sum} again to get rid of the summation over $i'$, we get
\begin{eqnarray*}
W(n,i,j) &=&\sum_{p=0}^{n-2} \sum_{a=1}^{i-1} \sum_{b=a+1}^{i} (-1)^{p+j+n} \binom{p}{n-j} A_{n,a,b} \\ & &
\qquad\qquad \qquad\times \left( \binom{b+n-1-a-i}{n-a-p-1} - \binom{n-1-i}{n-a-p-1} \right) \\ & &
+ \sum_{p=0}^{n-2} \sum_{a=1}^{i-1} \sum_{b=i+1}^{n} (-1)^{p+j+n} \binom{p}{n-j} A_{n,a,b}  \\ & &
\qquad\qquad \qquad\times\left( \binom{i+n-a-b}{n-a-b-p+i} - \binom{n-b}{n-a-b-p+i} \right) \\ & &
+ \sum_{p=0}^{n-2} \sum_{b=2}^{i} \sum_{a=1}^{b-1} (-1)^{p+j+n} \binom{p}{n-j} \binom{n-1-i}{p-i+a} A_{n,a,b} \\ & &
+ \sum_{p=0}^{n-2} \sum_{b=i+1}^{n} \sum_{a=1}^{i} (-1)^{p+j+n} \binom{p}{n-j} \binom{n-b}{p-i+a} A_{n,a,b},
\end{eqnarray*}
which simplifies to 
\begin{eqnarray*}
W(n,i,j) &=& 
\sum_{p=0}^{n-2} \sum_{a=1}^{i-1} \sum_{b=a+1}^{i} (-1)^{p+j+n} \binom{p}{n-j} \binom{b+n-1-a-i}{b-i+p} A_{n,a,b} \\
& & + \sum_{p=0}^{n-2} \sum_{a=1}^{i} \sum_{b=i+1}^{n} (-1)^{p+j+n} \binom{p}{n-j} \binom{i+n-a-b}{p} A_{n,a,b}.
\end{eqnarray*}
Now use the summation formula 
$$
\sum_{k} \binom{r}{m+k} \binom{s+k}{n} (-1)^k = (-1)^{r+m} \binom{s-m}{n-r}
\qquad (r,m,n \textrm{ integers}, r \ge 0)
$$
(which follows from the Chu-Vandermonde summation together with the basic transformation formulas
$\binom{a}{b}=\binom{a}{a-b}=(-1)^b \binom{b-a-1}{b}$) to eliminate the summation over $p$ and get
that $W(n,i,j)$ is equal to
\begin{multline*}
\sum_{a=1}^{i-1} \sum_{b=a+1}^{i} (-1)^{a+j+1} \binom{i-b}{j-1-a} A_{n,a,b} 
+ \sum_{a=1}^{i} \sum_{b=i+1}^{n} (-1)^{i+j+a+b} \binom{0}{a+b-i-j} A_{n,a,b} \\
=\ \ \  \sum_{a=1}^{i-1} \sum_{b=a+1}^{i} (-1)^{a+j+1} \binom{i-b}{j-1-a} A_{n,a,b} + 
\sum_{a=\max(1,i+j-n)}^{\min(i,j-1)} A_{n,a,i+j-a}.
\qquad\quad
\end{multline*}
Finally, plugging this into \eqref{eq:Fexpressedas} we obtain a formula for the coefficients $F_{n,i,j}$.
\begin{eqnarray*}
\nonumber
\!F_{n,i,j} &\!\!\!\!=& \!\!\![i \not= n] \left( \sum_{a=1}^{i-1} \sum_{b=a+1}^{i} (-1)^{a+j+1} \binom{i-b}{j-1-a} A_{n,a,b} + 
\sum_{a=\max(1,i+j-n)}^{\min(i,j-1)} A_{n,a,i+j-a} \right) \\ \nonumber
& & \qquad - \sum_{a=1}^{i-2} \sum_{b=a+1}^{i-1} (-1)^{a+j} \binom{i-1-b}{j-a} A_{n,a,b} \\ \nonumber & & \qquad - 
\sum_{a=\max(1,i+j-n)}^{\min(i-1,j)} A_{n,a,i+j-a} + [i=n-1][j=1] A_{n-1},
\end{eqnarray*}
which we rewrite in the slightly more convenient form
\begin{eqnarray} \nonumber
F_{n,i,j} &\!\!\!\!=& [i < j] A_{n,i,j} - [j < i] A_{n,j,i} + [i=n-1][j=1] A_{n-1} \\ \nonumber
& & \qquad + [i \not= n] \left( \sum_{a=1}^{i-1} \sum_{b=a+1}^{i} (-1)^{a+j+1} \binom{i-b}{j-1-a} A_{n,a,b} \right) \\ & & \qquad + 
\sum_{a=1}^{i-2} \sum_{b=a+1}^{i-1} (-1)^{a+j+1} \binom{i-1-b}{j-a} A_{n,a,b}.
\label{ilse}
\end{eqnarray}
We have attained our first important goal, which was to derive a formula that expresses $F_{n,i,j}$ in terms of
the doubly-refined enumeration numbers $\{A_{n,p,q}\}_{p,q}$.
This formula immediately implies ``half'' of Theorem \ref{thm-coefficients}, namely that $F_{n,i,j}=A_{n,i,j}$ for
$ 1\le i<j\le n$: If $a < b \le i < j$ then $i-b < j-1-a$ and 
$i-1-b < j - a$ and, consequently, $\binom{i-b}{j-1-a} = 0$ and $\binom{i-1-b}{j-a}=0$ in this case. 
Moreover, 
$1 \le i < j \le n$ implies $j \ge 2$ and thus $[j=1]=0$.

We will later use \eqref{ilse} again to prove the near-symmetry property of the $\hat{A}_{n,i,j}$'s (Theorem \ref{thm-almostsym}).

\section{Finishing the proof of Theorem \ref{thm-coefficients}}

We turn now to proving that $F_{n,i,j}=\hat{A}_{n,i,j}$ for the other case when $1\le j\le i\le n$. Because of the way the extended doubly-refined enumeration numbers were defined in \eqref{eq:extendednumbers}, the proof of Theorem \ref{thm-coefficients} will be complete once we prove the following lemma.

\begin{lemma} \label{rep_cnij_lemma}
For each $1\le j\le i \le n$ we have
\begin{equation}\label{eq:rep_cnij}
 F_{n,i,j} = \sum_{0\le p < q \le n-1} \ga_{i,j,p,q} F_{n,p,q},
 \end{equation}
where $\ga_{i,j,p,q}$ are defined in \eqref{eq:gammas}.
\end{lemma}
For this, we appeal to another identity satisfied by $\alpha_n$. The following lemma is proved in \cite{fischer1} (see Corollary 1 in that paper and the paragraph below it).

\begin{lem}
For each $1\le i \neq j\le n$ let $S_{i,j}$ be the operator that swaps the variables $k_i$ and $k_j$, i.e.,
$$ (S_{i,j} f)(k_1,\ldots,k_n) = f(k_1,\ldots,k_{i-1},k_j,k_{i+1},\ldots,k_{j-1},k_i,k_{j+1},\ldots,k_n) \quad (\text{when }i<j).$$
Then for each $1\le i\le n-1$, $\alpha_n$ satisfies the identity
$$ (\text{Id}+S_{i,i+1})(\text{Id}+E_{k_i} E_{k_{i+1}} - E_{k_{i+1}}) \alpha_n \equiv 0 $$
(where $\text{Id}$ is the identity operator). More explicitly this can be written as the 6-term identity
\begin{multline}
\alpha_n(\ldots k_i,k_{i+1}\ldots)+\alpha_n(\ldots k_i+1,k_{i+1}+1\ldots)-\alpha_n(\ldots k_i,k_{i+1}+1\ldots) \\ =
-\alpha_n(\ldots k_{i+1},k_i\ldots) - \alpha_n(\ldots k_{i+1}+1,k_i+1\ldots) + \alpha_n(\ldots k_{i+1},k_i+1\ldots).
\label{eq:familyidentity1}
\end{multline}
\end{lem}

Applying \eqref{eq:familyidentity1} with $i=n-1$ gives the following lemma.

\begin{lemma} The function $G_n(x,y)$ defined in \eqref{eq:definitionG} satisfies
\begin{multline}\label{eq:secondsystem}
G_n(x,y)+G_n(x+1,y+1)-G_n(x,y+1) \\
= -G_n(y+1,x-1) - G_n(y+2,x) + G_n(y+1,x).
\end{multline}
\end{lemma}

\begin{lemma} For each $n\ge1$, the coefficients $F_{n,i,j}$ satisfy the following system of linear equations.
\begin{equation}\label{eq:lineareqns_secondsystem}
F_{n,i,j} + \sum_{p = i+1}^{n} \sum_{q=j}^{n} F_{n,p,q}
= -F_{n,j,i} - \sum_{p = i}^{n} \sum_{q=j+1}^{n} F_{n,q,p} \qquad (1\le i,j\le n).
\end{equation}
\end{lemma}

\begin{proof} 
The claim will follow by expanding all 6 terms of eq. \eqref{eq:secondsystem} in the linear basis
$ \left\{ \binom{x+i-1}{i-1}\binom{y+j}{j-1}\right\}_{i,j}$, using \eqref{sum}.
For the terms on the left-hand side, we have that
\begin{eqnarray*}
G_n(x,y)&=& \sum_{i,j} F_{n,i,j} \binom{x+i-1}{i-1}\binom{y+j}{j-1}, \\
G_n(x,y+1) &=& \sum_{i,q} F_{n,i,q} \binom{x+i-1}{i-1}\binom{y+q+1}{q-1}
\\ &=& \sum_{i,q} F_{n,i,q} \binom{x+i-1}{i-1} \left(\sum_{j=1}^q \binom{y+j}{j-1} \right)
\\ &=& \sum_{i,j} \left(\sum_{q=j}^n F_{n,i,q}\right) \binom{x+i-1}{i-1} \binom{y+j}{j-1},
%\\
\end{eqnarray*} \begin{eqnarray*}
G_n(x+1,y+1) &=& \sum_{p,q} F_{n,p,q} \binom{x+p}{p-1}\binom{y+q+1}{q-1}
\\ &=& \sum_{p,q} F_{n,p,q} \left(\sum_{i=1}^p \binom{x+i-1}{i-1}\right) \left(\sum_{j=1}^q \binom{y+j}{j-1} \right)
\\ &=& \sum_{i,j} \left(\sum_{p=i}^n \sum_{q=j}^n F_{n,p,q} \right) \binom{x+i-1}{i-1}\binom{y+j}{j-1},
\end{eqnarray*}
and similarly for the terms on the right-hand side we have
\begin{eqnarray*}
G_n(y+1,x-1) &=& \sum_{i,j} F_{n,j,i} \binom{x+i-1}{i-1}\binom{y+j}{j},
\\ G_n(y+1,x) &=& \sum_{i,j} \left(\sum_{p=i}^n F_{n,j,p} \right) \binom{x+i-1}{i-1}\binom{y+j}{j},
\\ G_n(y+2,x) &=& \sum_{i,j} \left(\sum_{p=i}^n \sum_{q=j}^n F_{n,q,p} \right) \binom{x+i-1}{i-1}\binom{y+j}{j}.
\end{eqnarray*}
By combining these expansions with \eqref{eq:secondsystem} and comparing coefficients on both sides we get \eqref{eq:lineareqns_secondsystem}.
\end{proof}

\begin{proof}[Proof of Lemma \ref{rep_cnij_lemma}]
Consider the system \eqref{eq:lineareqns_secondsystem} as an inhomogeneous system of equations in indeterminate variables $(F_{n,i,j})_{1\le j\le i\le n}$, where $(F_{n,i,j})_{1\le i< j\le n}$ are given and considered as constants. We want to show that this system has a unique solution that is given by \eqref{eq:rep_cnij}. Note that the equations \eqref{eq:lineareqns_secondsystem} are symmetric in $i$ and $j$, so one may consider only the equations indexed by pairs $(i,j)$ with $j\le i$. First, we transform the system of equations into an equivalent one that is more convenient to handle: For each $1\le i\le n-1$, subtract the equation with index $(i+1,i+1)$ from the equation with index $(i,i)$, and for each $1\le j<i\le n$, subtract the equation with index $(i,j+1)$ from the equation with index $(i,j)$. We obtain the following equivalent system:
\begin{eqnarray}
F_{n,n,n} &=& 0, \label{eq:systemeq1} \\
\sum_{k=i}^n F_{n,k,i} + \sum_{k=i+2}^n F_{n,i+1,k} &=& 0 \qquad (1\le i\le n-1), \label{eq:systemeq2} \\
\sum_{k=i}^n F_{n,k,j} - F_{n,i,j+1} + F_{n,j,i} + \sum_{k=i+1}^n F_{n,j+1,k} &=&0 \qquad (1\le j<i\le n).
\label{eq:systemeq3}
\end{eqnarray}
It is now easy to see that this system of equations is triangular (with all coefficients equal to 1 on the diagonal) when the variables are ordered in a reverse lexicographical order, scanning the rows from 
bottom to top ($i=n$ to $i=1$) and each row from right to left ($j=i$ to $j=1$). This implies that the system has a unique solution. It remains to verify that \eqref{eq:rep_cnij} is in fact a solution. Eq. \eqref{eq:systemeq1} is trivial (and actually can be thought of as the case $i=n$ of \eqref{eq:systemeq2}), so we verify \eqref{eq:systemeq2} and \eqref{eq:systemeq3}.

For  \eqref{eq:systemeq2}, substitute \eqref{eq:rep_cnij} into the equation and equate the coefficients of $F_{n,p,q}$ to 0 for each $p<q$, to see that it is necessary to check that
\begin{equation}
\sum_{k=i}^n c_{k,i,p,q} = -[p=i+1] \qquad (1\le i\le n-1, \ 1\le p<q\le n).
\label{eq:needtocheck}
\end{equation}
Note that $c_{i,j,p,q}$ can be written as 
$$ c_{i,j,p,q} = (-1)^{i+q+1}\left( \binom{p-j+1}{q-i}[p\ge j] - \binom{p-j-1}{q-i-1}[p\ge j+1] \right), $$
so to check \eqref{eq:needtocheck}, we divide into 4 cases:
\begin{enumerate}
\item $p<i$: In this case we get immediately that
$ \sum_{k=i}^n c_{k,i,p,q} = \sum_{k=i}^n 0 = 0. $
\item $p=i$: 
\begin{eqnarray*}
\sum_{k=i}^n c_{k,i,p,q} &=& \sum_{k=i}^n (-1)^{k+q+1} \binom{1}{q-k} \\ &=& (-1)^{q+1} \left((-1)^{q-1}\binom{1}{1}+(-1)^q\binom{1}{0}
\right) = 0.
\end{eqnarray*}
\item $p=i+1$:
\begin{eqnarray*}
\sum_{k=i}^n c_{k,i,p,q} &=&
\sum_{k=i}^n (-1)^{k+q+1} \left( \binom{2}{q-k}-\binom{0}{q-k-1}\right) \\ &=&
-\binom{2}{2}+\binom{2}{1}-\binom{2}{0}-\binom{0}{0} = -1 = -[p=i+1]. 
\end{eqnarray*}
\item $p\ge i+2$:
\begin{eqnarray*}
\sum_{k=i}^n c_{k,i,p,q} &=&
\sum_{k=i}^n (-1)^{q+k+1}\left( \binom{p-i+1}{q-k}-\binom{p-i-1}{q-k-1}\right) \\ &=&
- \sum_{m=0}^{p-i+1}(-1)^m \binom{p-i+1}{m} + \sum_{m=0}^{p-i-1} (-1)^m \binom{p-i-1}{m}  = 0.
\end{eqnarray*}
\end{enumerate}
This confirms \eqref{eq:systemeq2}. The verification of \eqref{eq:systemeq3} is based on a similar case analysis, where after substituting \eqref{eq:rep_cnij} into \eqref{eq:systemeq3} one has to verify the identity
$$ -c_{i,j+1,p,q} + \sum_{k=i}^n c_{k,j,p,q} = -[p=j][q=i] - [p=j+1][q\ge i+1]
$$
for $1 \le j < i \le n$ and $1\le p<q\le n$. Here are the different cases that require checking.
\begin{enumerate}
\item $(p,q)=(j,i)$
\item $(p,q)=(j+1,i+1)$
\item $p=j+1,q\ge i+2$
\item $p=j, q<i$
\item $p=j, q>i$
\item $p=j+1, q\le i$
\item $p<j$
\item $p\ge j+2$
\end{enumerate}
We omit the details of this verification, which are straightforward and easy to fill in.
\end{proof}

With Lemma \ref{rep_cnij_lemma} proved, the proof of Theorem \ref{thm-coefficients} is also complete.

\section{Proofs of Theorem \ref{mainresult} and Proposition \ref{two-conjectures}\label{sectionproofmain}}

We have identified the coefficients $F_{n,i,j}$ for $1\le i<j\le n$ as the extended doubly-refined enumeration numbers $\hat{A}_{n,i,j}$. 
To prove Theorem \ref{mainresult}, we now derive the system of linear equations satisfied by the $F_{n,i,j}$'s. This will follow from several of the symmetry properties of $\alpha_n$. Two of them are the combinatorially obvious identities:
\begin{eqnarray}
\alpha_n(k_1 + t, k_2 + t, \ldots, k_n + t) &=& \alpha_n(k_1, k_2, \ldots, k_n), \label{eq:alphanidentity1} \\
\alpha_n(-k_n,\ldots,-k_2,-k_1) &=& \alpha_n(k_1,k_2,\ldots,k_n). \label{eq:alphanidentity2} 
\end{eqnarray}
In \cite[Lemma 5]{fischer2} the following additional identity is proved.

\begin{lemma}
\begin{equation}
\alpha_n(k_2,k_3,\ldots,k_n,k_1-n)= (-1)^{n-1} \alpha_n(k_1,k_2,\ldots,k_n). \label{eq:alphanidentity3}
\end{equation}
\end{lemma}

As a consequence of the three identities \eqref{eq:alphanidentity1}, \eqref{eq:alphanidentity2}, \eqref{eq:alphanidentity3}, we get:
\begin{lemma} \label{polynomialidentity} The polynomial $G_n(x,y)$ from \eqref{eq:definitionG} satisfies
\begin{equation}\label{eq:poly_identity1}
G_n(x,y) = G_n(-2n-y, -2n-x).
\end{equation}
\end{lemma}

\begin{proof} We use the identities \eqref{eq:alphanidentity1}, \eqref{eq:alphanidentity2}, and \eqref{eq:alphanidentity3} to get
\begin{eqnarray*}
G_n(x,y)&=& \alpha_n(1,\ldots,n-2,n-1+x,n+y) \\ &=& \alpha_n(-n-y,-n+1-x,-n+2,-n+3,\ldots,-2,-1)
\\ &=& (-1)^{n-1} \alpha_n(-n+1-x,-n+2,\ldots,-1,-2n-y)
\\ &=& (-1)^{n-1} (-1)^{n-1} \alpha_n(-n+2,\ldots,-2,-1,-2n-y,-2n+1-x)
\\ &=& \alpha_n(1,2,3,\ldots,n-2,-n-1-y,-n-x) = G_n(-2n-y,-2n-x).
\end{eqnarray*}
\end{proof}

\begin{proof}[Proof of Theorem \ref{mainresult}]
Expand both sides of \eqref{eq:poly_identity1} using \eqref{eq:uniqueexp} and Theorem \ref{thm-coefficients}, and use the Chu-Vandermonde identity
to get that
\begin{samepage}
\begin{eqnarray*} \sum_{i,j} \hat{A}_{n,i,j} \binom{x+i-1}{i-1}\binom{y+j}{j-1}
&=& \sum_{p,q}  \hat{A}_{n,p,q} \binom{-2n-y+p-1}{p-1}\binom{-2n-x+q}{q-1} \qquad \qquad\qquad
\end{eqnarray*}

\vspace{-12.0pt}
\begin{eqnarray*}
%\qquad \qquad 
%\qquad 
&=& \sum_{p,q} \hat{A}_{n,q,p} \binom{-2n-x+p}{p-1}\binom{-2n-y+q-1}{q-1}
\\ &=& \sum_{p,q} \hat{A}_{n,q,p} \left( \sum_{i=1}^p (-1)^{i-1} \binom{p-2n+1}{p-i}\binom{x+i-1}{i-1}\right)
\\ & & \qquad \qquad \times\left( \sum_{j=1}^q (-1)^{j-1} \binom{q-2n+1}{q-j}\binom{y+j}{j-1} \right)
\\ &=& \sum_{i,j=0}^{n-1} (-1)^{i+j} \left( \sum_{p=i}^n \sum_{q=j}^n \hat{A}_{n,q,p} \binom{p-2n+1}{p-i}\binom{q-2n+1}{q-j} \right)
\binom{x+i-1}{i-1} \binom{y+j}{j-1}
\\ &=& \sum_{i,j=0}^{n-1} \left( \sum_{p=i}^n \sum_{q=j}^n (-1)^{p+q}  \binom{2n-i-2}{p-i}\binom{2n-j-2}{q-j} \hat{A}_{n,q,p} \right)
\binom{x+i-1}{i-1} \binom{y+j}{j-1}.
\end{eqnarray*}
\end{samepage}

\noindent
Comparing coefficients of $\binom{x+i-1}{i-1}\binom{y+j}{j-1}$ gives exactly the equations \eqref{eq:lineareqns}.
\end{proof}

\begin{proof}[Proof of Proposition \ref{two-conjectures}]
The technique used to prove Lemma \ref{polynomialidentity} generalizes easily to the polynomial $G_n(x_1,\ldots,x_d)$, and gives
the identity
$$ G_n(x_1,\ldots,x_d) = (-1)^{(n-1)d} G_d(-2n-x_d,\ldots,-2n-x_2,-2n-x_1). $$
Assuming Conjecture \ref{conj-drefined}, the linear equations \eqref{eq:drefinedeqs} follow by expanding the functions on both sides of this equation in the basis
$$ \left\{ 
\binom{x_1+j_1-1}{j_1-1}
\binom{x_2+j_2}{j_2-1}
\ldots
\binom{x_d+j_d+d-2}{j_d-1}
\right\}_{j_1,j_2,\ldots,j_d\ge 1}
$$
and equating coefficients. The details are similar to the computation in the proof above and are omitted.
\end{proof}

\section{Proof of Theorem \ref{thm-almostsym}}

First, note that Theorem \ref{thm-almostsym} is equivalent to the evaulation $\hat{A}_{n,n,2}=A_{n-2}-A_{n-1}$ together with the claim that the modified numbers defined by 
$$\tilde{A}_{n,i,j}=\hat{A}_{n,i,j}-[i=n-1][j=1]A_{n-1}$$ 
satisfy the symmetry 
\begin{equation}\label{eq:modifiedsymmetry}
\tilde{A}_{n,i,j}=\tilde{A}_{n,n+1-j,n+1-i}
\end{equation}
for \dd{all} $1\le i,j\le n$. The first claim is a special case of \eqref{eq:interestingid} (which was a direct consequence of  \eqref{eq:extendednumbers}). It remains to prove the symmetry. 
We already know it for $i<j$, so we show this for  $i\ge j$.

First,
consider the case $i \not= n$ and $j \not=1$. In this case, by \eqref{ilse} we have 
\begin{eqnarray*}
\hat{A}_{n,i,j}&=& [i < j] A_{n,i,j} - [j < i] A_{n,j,i}  \\
& & + \sum_{a=1}^{i-1} \sum_{b=a+1}^{i} (-1)^{a+j+1} \binom{i-b}{j-1-a} A_{n,a,b} \\ & & + 
\sum_{a=1}^{i-2} \sum_{b=a+1}^{i-1} (-1)^{a+j+1} \binom{i-1-b}{j-a} A_{n,a,b}.
\end{eqnarray*}
Moreover, by \eqref{eq:extendednumbers} and since $\binom{n}{k}_{+} = 0$ for
$n < k$, if $i\ge j$ (including if $i=n$ or $j=1$) we have 
\begin{eqnarray}
\hat{A}_{n,i,j}  &\!\!=\!\!& \sum_{j \le a < b \le \min(i-j+1+a,n)}
(-1)^{i+b+1} \left( \binom{a-j+1}{b-i}_{+} - \binom{a-j-1}{b-i-1}_{+} \right)
A_{n,a,b} \nonumber \\ 
 &\!\!=\!\!& A_{n,j,i+1} + \sum_{a=j}^{n} \sum_{b=a+1}^{\min(i-j+1+a,n)} \!\!\!\! (-1)^{i+b+1} 
\left( \binom{a-j+1}{b-i} - \binom{a-j-1}{b-i-1} \right) A_{n,a,b} \nonumber \\
&\!\!=\!\!& A_{n,j,i+1} + \sum_{a=j}^{n} \sum_{b=a+1}^{\min(i-j+1+a,n)} \!\!\!\! (-1)^{i+b+1} 
\left( \binom{a-j-1}{b-i} + \binom{a-j}{b-i-1} \right) A_{n,a,b} \nonumber \\
&\!\!=\!\!&  - [i \not= j] A_{n,j,i} 
 + \sum_{a=j+1}^{n} \sum_{b=a+1}^{n} (-1)^{i+b+1} \binom{a-j-1}{b-i} A_{n,a,b} \nonumber \\
& & \qquad  + \sum_{a=j}^{n} \sum_{b=a+1}^{n} (-1)^{i+b+1} \binom{a-j}{b-i-1} A_{n,a,b}.   
\label{dan}
\end{eqnarray}
If we add these two representations for $F_{n,i,j}$ we get that 
\begin{eqnarray}
2\hat{A}_{n,i,j} &=& - 2[i \not= j] A_{n,j,i} 
+ \sum_{a=j+1}^{n} \sum_{b=a+1}^{n} (-1)^{i+b+1} \binom{a-j-1}{b-i} A_{n,a,b} 
\nonumber \\ & &  
 + \sum_{a=j}^{n} \sum_{b=a+1}^{n} (-1)^{i+b+1} \binom{a-j}{b-i-1} A_{n,a,b} \nonumber \\
& &  + \sum_{a=1}^{i-1} \sum_{b=a+1}^{i} (-1)^{a+j+1} \binom{i-b}{j-1-a} A_{n,a,b} \nonumber \\ & &+ 
\sum_{a=1}^{i-2} \sum_{b=a+1}^{i-1} (-1)^{a+j+1} \binom{i-1-b}{j-a} A_{n,a,b}. 
\label{eq:becauseof}
\end{eqnarray}
We set 
\begin{eqnarray*}
D_{n,i,j} &=& \sum_{a=j+1}^{n} \sum_{b=a+1}^{n} (-1)^{i+b+1} \binom{a-j-1}{b-i} A_{n,a,b} 
\\ & & + \sum_{a=j}^{n} \sum_{b=a+1}^{n} (-1)^{i+b+1} \binom{a-j}{b-i-1} A_{n,a,b} \\
& &  + \sum_{a=1}^{i-1} \sum_{b=a+1}^{i} (-1)^{a+j+1} \binom{i-b}{j-1-a} A_{n,a,b} \\ & &+ 
\sum_{a=1}^{i-2} \sum_{b=a+1}^{i-1} (-1)^{a+j+1} \binom{i-1-b}{j-a} A_{n,a,b}
\end{eqnarray*}
and observe that because of \eqref{eq:becauseof}, it is enough to show that $D_{n,i,j}=D_{n,n+1-j,n+1-i}$.
In this formula, we replace all occurrences of $A_{n,a,b}$ by $A_{n,n+1-b,n+1-a}$.
Then we replace $b$ by $n+1-a'$ and $a$ by $n+1-b'$ and obtain the following.
\begin{eqnarray*}
D_{n,i,j} &=& \sum_{b'=1}^{n-j} \sum_{a'=1}^{b'-1} (-1)^{i+a'+n} \binom{n-b'-j}{n+1-i-a'} A_{n,a',b'} \\
& & + \sum_{b'=1}^{n+1-j} \sum_{a'=1}^{b'-1} (-1)^{i+a'+n} \binom{n+1-j-b}{n-i-a} A_{n,a',b'} \\
& &+ \sum_{b'=n+2-i}^{n} \sum_{a'=n+1-i}^{b'-1} (-1)^{j+b'+n} \binom{-n-1+i+a}{-n-2+j+b'} A_{n,a',b'} \\
& & + \sum_{b'=n+3-i}^{n} \sum_{a'=n+2-i}^{b'-1} (-1)^{j+b'+n} \binom{-n-2+i+a'}{-n-1+j+b'} A_{n,a',b'}
\end{eqnarray*} 
If we exchange the order of the summation in the four double sums (and replace $a'$ by $a$ and $b'$ by $b$)
then we obtain 
\begin{eqnarray*}
D_{n,i,j} &=&  \sum_{a=1}^{n-j-1} \sum_{b=a+1}^{n-j} (-1)^{i+a+n} \binom{n-b-j}{n+1-i-a} A_{n,a,b} \\
& & + \sum_{a=1}^{n-j} \sum_{b=a+1}^{n+1-j} (-1)^{i+a+n} \binom{n+1-j-b}{n-i-a} A_{n,a,b} \\
& & + \sum_{a=n+1-i}^{n-1} \sum_{b=a+1}^{n} (-1)^{j+b+n} \binom{-n-1+i+a}{-n-2+j+b} A_{n,a,b} \\
& & + \sum_{a=n+2-i}^{n-1} \sum_{b=a+1}^{n} (-1)^{j+b+n} \binom{-n-2+i+a}{-n-1+j+b} A_{n,a,b}
\end{eqnarray*} 
and this expression is obviously equal to $D_{n,n+1-j,n+1-i}$.

\medskip

Finally, we prove \eqref{eq:modifiedsymmetry} 
for the case that $i=n$ or $j=1$. It suffices to consider the 
case $i=n$, since $j=1$ implies $n+1-j=n$. That is, we have to show 
$$
\hat{A}_{n,n,j} - \hat{A}_{n,n+1-j,1} + [j=2] A_{n-1} = 0.
$$
For $\hat{A}_{n,n,j}$ we use the formula from \eqref{dan} and for $\hat{A}_{n,n+1-j,1}$ we use \eqref{ilse}. This gives
\begin{samepage}
\begin{equation*}
\hat{A}_{n,n,j} - \hat{A}_{n,n+1-j,1} + [j=2] A_{n-1}
\qquad\qquad\qquad\qquad\qquad\qquad\qquad\qquad\qquad\qquad\qquad
\end{equation*}

\vspace{-20.0pt}
\begin{eqnarray*}
&=&
-[n \not= j] A_{n,j,n} + \sum_{a=j+1}^{n} \sum_{b=a+1}^{n} (-1)^{n+b+1} \binom{a-j-1}{b-n} A_{n,a,b} \\
& & + \sum_{a=j}^{n} \sum_{b=a+1}^{n} (-1)^{n+b+1} \binom{a-j}{b-n-1} A_{n,a,b} 
\\ & & + [n+1-j \not= 1] A_{n,1,n+1-j} 
- [j=2] A_{n-1} + [j=2] A_{n-1}\\
& & - \sum_{a=1}^{n-j} \sum_{b=a+1}^{n+1-j} (-1)^a \binom{n+1-j-b}{-a} A_{n,a,b} 
-\!\! \sum_{a=1}^{n-1-j} \sum_{b=a+1}^{n-j} (-1)^a \binom{n-j-b}{1-a} A_{n,a,b}.
% \\& & 
\end{eqnarray*}
\end{samepage}
Since the second and the third double sums vanish and $[n \not= j] A_{n,j,n} = [n+1-j \not= 1] A_{n,1,n+1-j} $, it remains 
to show that
$$
- \sum_{a=j+1}^{n-1} A_{n,a,n} + \sum_{b=2}^{n-j} A_{n,1,b} = 0.
$$
This follows from the symmetry of the $A_{n,p,q}$'s when $p<q$.
\qed

\section*{Appendix A: Numerical tables\label{numerical}}

\begin{center}
\begin{figure}[h!]
\hspace{50.0pt}
\begin{tabular}{*{13}{c@{\hspace{0.02in}}}}
&&&&&& 1 &&&&&& \\
&&&&& 1 && 1 &&&&& \\
&&&& 2 && 3 && 2 &&&& \\
&&& 7 && 14 && 14 && 7 &&& \\
&& 42 && 105 && 135 && 105 && 42 && \\
& 429 && 1287 && 2002 && 2002 && 1287 && 429 & \\
7436 && 26026 && 47320 && 56784 && 47320 && 26026 && 7436 
\end{tabular}
\caption{The numbers $A_{n,k}$ for $1\le n \le 7$}
\end{figure}
\end{center}

\begin{center}
\begin{figure}[h!]
$$
\begin{array}{ll}
\begin{array}{l}n=3 \\ \ \\ \ \end{array}
& \left( \begin{array}{ccc}0&1&1\\1&1&1\\-2&-1&0\end{array}\right) \\ \ \\
\begin{array}{l}n=4 \\ \ \\ \  \\ \ \\\end{array}
 & \left( \begin{array}{cccc}0&2&3&2\\-2&2&4&3\\2&1&2&2\\-7&-5&-2&0
\end{array}\right) \\ \ \\
\begin{array}{l}n=5 \\ \ \\ \  \\ \ \\ \ \\  \end{array}
 & \left( \begin{array}{ccccc} 0&7&14&14&7\\-7&7&23&26&14\\
 -21&-2&16&23&14\\
 7&-7&-2&7&7\\
 -42&-35&-21&-7&0
\end{array}\right) \\ \ \\
\begin{array}{l}n=6 \\ \ \\ \  \\ \ \\ \ \\ \ \\  \end{array}
 & \left( \begin{array}{cccccc} 
 0&42&105&135&105&42\\
 -42&42&203&300&250&105\\
 -147&-56&161&322&300&135\\
 -282&-179&-8&161&203&105\\
 42&-177&-179&-56&42&42\\
 -429&-387&-282&-147&-42&0
\end{array}\right) \\ \ \\
\begin{array}{l}n=7 \\ \ \\ \  \\ \ \\ \ \\ \ \\ \ \\\end{array}
 & \left( \begin{array}{cccccccc} 
 0 & 429 & 1287 & 2002 & 2002 & 1287 & 429 \\
 -429 & 429 & 2847 & 5174 & 5551 & 3731 & 1287 \\
 -1716 & -1131 & 2418 & 6422 & 7748 & 5551 & 2002 \\
 -3718 & -3874 & -546 & 4004 & 6422 & 5174 & 2002 \\
 -5720 & -5707 & -4043 & -546 & 2418 & 2847 & 1287 \\
 429 & -4433 & -5707 & -3874 & -1131 & 429 & 429 \\
 -7436 & -7007 & -5720 & -3718 & -1716 & -429 & 0
\end{array}\right) \\
\end{array}
$$
\caption{The extended doubly-refined enumeration coefficient matrices $(\hat{A}_{n,i,j})_{i,j=1}^n$ for $3\le n\le 7$.}
\end{figure}
\end{center}

\newpage

\ 

\newpage

\section*{Appendix B: The \texttt{Mathematica} package \texttt{RefinedASM}}

The \texttt{Mathematica 6.0} package \texttt{RefinedASM} can be downloaded from the 
authors' web pages.
Here's a sample output to demonstrate the computation of the $\hat{A}_{n,i,j}$'s for $n=6$ 
(which at the same time verifies the correctness of Conjecture \ref{conj-sufficiency} for that value of $n$)
and the verification of Conjecture \ref{conj-explicitformula} for $3\le n\le 8$:

\vspace{10.0pt}
\includegraphics{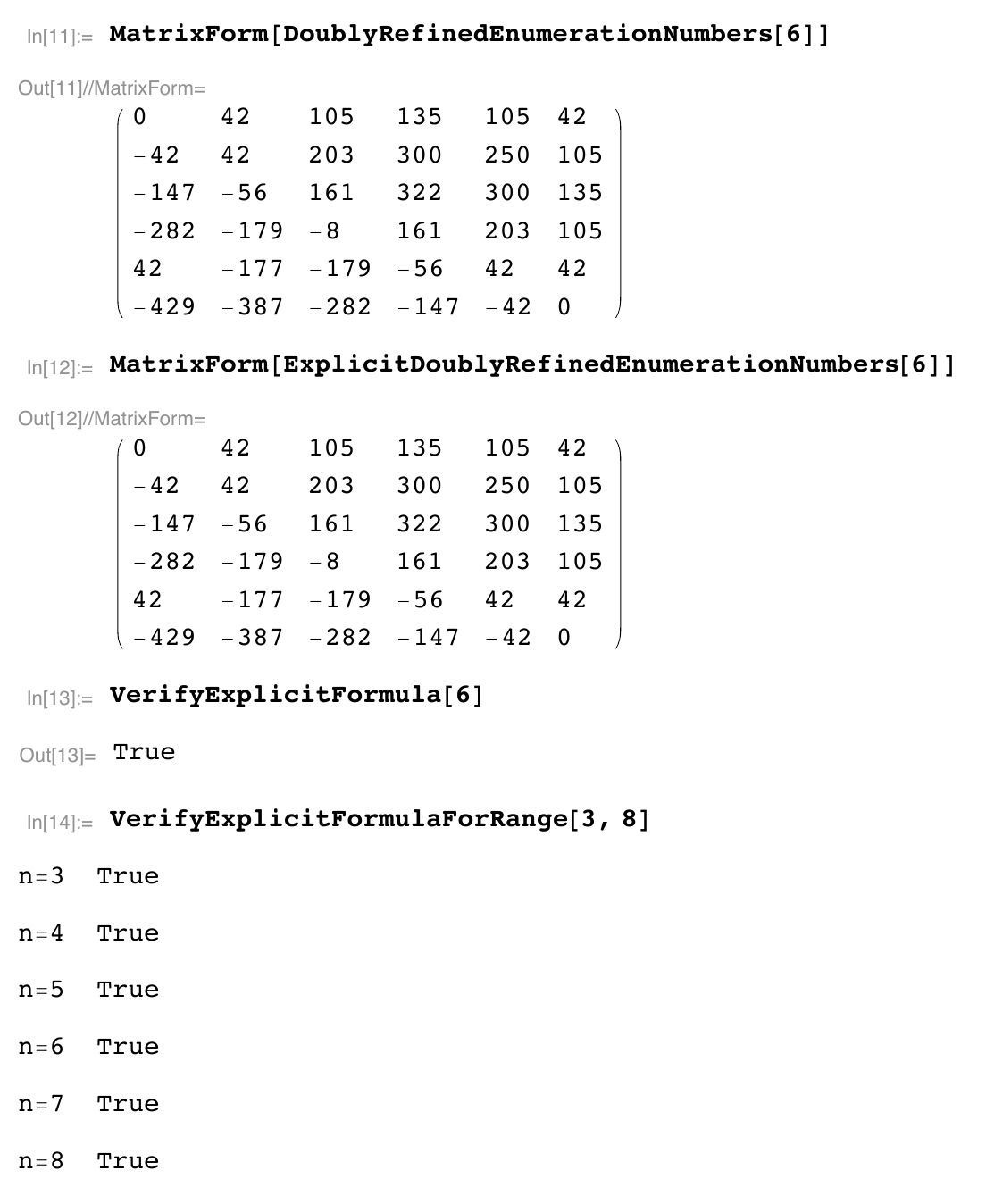}

\newpage
 \noindent
\textsc{Ilse Fischer} \\
\textsc{
Institut f\"ur Mathematik, Universit\"at Klagenfurt \\
9020 Klagenfurt, and \\
Fakult\"at f\"ur Mathematik, Universit\"at Wien \\
1090 Wien, Austria \\
}
\texttt{ilse.fischer@univie.ac.at}

\bigskip \noindent
\textsc{Dan Romik} \\
\textsc{
Einstein Institute of Mathematics, The Hebrew University \\
Givat-Ram, Jerusalem 91904, Israel \\
}
\texttt{romik@math.huji.ac.il}

\end{document}